\documentclass[11pt]{article}
\usepackage{amsmath,amssymb,amsthm,mathrsfs}

\begin{document}

\renewcommand{\baselinestretch}{1.3}
\renewcommand{\arraystretch}{1.3}

\begin{center}{\bf\LARGE  Frobenius Theorem  and a Principle for Critical Point under the Constraint of $C^1$ Submanifold in Banach Spaces }\\
\vskip 0.5cm Ma Jipu$^{1,2}$
\end{center}

{\bf Abstract}\quad   Let $\Lambda$ be an open set in a Banach space $E$ and $x_0$ be a point in $\Lambda$.
  We consider the family of subspaces in $E$,   $\mathcal{F}$$=\{M(x)\}_{x\in\Lambda}$,  especially, where  the
  dimension of $M(x)$ may be infinity,   and investigate the necessary and sufficient
condition for $\cal{F}$ being $c^1$  integrable at  $x_0$. For this purpose
we introduce the concept of the  co-final set $J(x_0,E_*)$ of
$\mathcal{F}$ at $x_0$, so that $M(x)$ for $x\in J(x_0,E_*)$ possesses an unique operator valued coordinate, $\alpha (x)$ $\in B(M(x_0), E_*)$. Then a Frobenius theorem in Banach
spaces is established.The theorem expands the well known Frobenius theorem with $ \dim M(x)\leq const.< \infty$ to $ \dim M(x)\leq \infty .$  When
$J(x_0,E_*)$ is trivial, the conditions of the theorem are reduced to the
solvability of an initial value problem for a differential equation
in Banach spaces. Specially let $f$ be $c^1$ map from an open set $U\subset E$
into another Banach space $F,N_x=N(f'(x)),$ and
$\mathcal{F}=\{N_x:\forall x\in U\}$. We prove the following
theorem: if $x_0\in U$  is a generalized
 regular point of $f$, then
the co-final set of $\mathcal{F}$ at $x_0$ is trivial, and
$\mathcal{F}$ is $c^1$ integrable at $x_0$. Owing to the generalized regular point being often the case in nonlinear functional analysis, the theorem presents
 several families of subspaces $\mathcal{F}=\{N_x:\forall x\in U\}$ with $\dim N_x\leq {\infty} $ and the trivial co-final set, all of which are $c^1$  integrable at $x_0$; while it provides several kinds of solvable differential equations with an initial value. Let $\Lambda =B(E,F)$ $\backslash $ $ \{0\},$ $ M(X)=\{T \in B(E,F):TN(X)\subset R(X)\}$ for $X\in \Lambda ,$  and the family of subspaces $\mathcal{F} = \{M(X):X\in \Lambda \}.$ Where $B(E,F)$ denotes all of linear bounded operators from $E$ into $F.$ We prove the theorem as follows, if $A\in \Lambda $ is double splitting, then $\mathcal{F} $ is a smooth integrable at $A$. Let $\Phi $ be any one of  $F_k, \Phi_{m,n}, \Phi_{m,\infty}$ and $\Phi_{\infty,n}$ ($m,n<\infty$). We have the following global results: $\Phi $ is a smooth submanifold in $B(E,F)$ and tangent to $M(X)$ at any $X\in \Phi .$ Especially where the co-final set of $\mathcal {F} $ at $A$ is non-trivial in generalized . These seem
to be useful for developing of differential topology, global
analysis and the geometrical method in differential equations.
Finally let $f:U \subset E\rightarrow (-\infty , \infty )$ be a $c^1$ nonlinear functional, and $S$ a $c^1$
 submanifold in $U.$ We give the result  as follows, if $x_0 \in S$ is a critical point of $f$ under the constraint $S$,
 then $N(f'(x_0))\supset T_{x_0}S.$ In view of the above Banach submanifolds with the explicit expression of the tangent space,the principle should be potential in application.
 \vskip 0.2cm
\begin{center}{\bf 1\quad   Introduction and Preliminary} \end{center}\vskip 0.2cm

Let $E$ be a Banach space and $\Lambda$ an open set in $E$. Assign a
subspace $M(x)$ in $E$ for each $x\in\Lambda$ where especially
$\dim M(x)$ may be infinity. We consider the family of subspaces
$\mathcal{F}=\{M(x):\forall x\in\Lambda\}$ and investigate the
necessary and sufficient condition for $\mathcal{F}$ being $c^1$
integrable at a point $x_0\in\Lambda$. When $E={R}^n$ or ${C}^n$,
the Frobenius theorem together with the implicit function theore
and the existence theorem for  solution of ordinate differential
equations are the three main pillars supporting differential
topology  and calculus on manifolds (refer [Abr]). Hence we have to
investigate Frobenius theorem in Banach spaces in which $\dim M(x)$
may be infinity.

For this purpose we introduce the concept of co-final set
$J(x_0,E_*)$ of $\mathcal{F}$ at $x_0$, so that it yields the following results:$M(x)$ for $x\in J(x_0,E_*)$ possesses an unique operator valued coordinate $\alpha \in B(M(x_0),E_*)$ (see Theorem 1.5); when $\mathcal {F}$ is $ c^1$ integrable at $x_0$, its integral submainfold $S$ satisfies $J(x_0,E_*) \supset S\bigcap U$ for some neighborhood $U_0$ at $x_0$ (see Theorem 2.1). Then a Frobenius theorem in Banach spaces is established (see Theorem 2.3). The theorem expands the well known Frobenius theorem with $\dim M(x) \leq const. < \infty $ to $\dim M(x)\leq \infty .$  When the co-final set of
$\mathcal{F}$ at $x_0$ is trivial (see Definition 3.1 in Section 3),
the conditions of the theorem are reduced to the solvability of an
initial value problem for a differential equation in Banach spaces.
Specially let $f$ be a $c^1$ map from an open set $U\subset E$ into another
Banach space $F, N_x=N(f'(x))$ for $x\in U$, and
$\mathcal {F}=\{N_x:\forall x\in U\}$. Let us make the following appointments with $c^1$ map $f:$ $x_0$ is said to be a Fredholm point of $f$ provided $\dim N( f'(x))=const.<\infty $ and $\mathrm{codim} R(f'(x) )=const.< \infty$ for $x$ near $x_0 ;$ $x_0$ is said to be a semi-Fredholm point of $f$ provided $\dim N(f'(x))= const.<\infty$ and $\mathrm{codim} R(f'(x))=\infty,$ or $\dim N(f'(x))=\infty $ and $\mathrm{codim} R(f'(x))=const.< \infty $ for $x$ near $x_0.$ Recall that $x_0\in U$ is said to be a generalized regular point of $f$ provided $x_0$ is a locally fine point of $f'(x).$  We  prove the theorem as follows, if $x_0\in U$
is a generalized regular point, then
the co-final set of $\mathcal{F}$ at $x_0$ is trivial, and $\mathcal {F}$ is $c^1$ integrable at $x_0$ (see Theorem 3.2). By Theorem 1.2 below it is easy to observe that all the following points of ${f}$ are generalized regular point of $f$: immersion, subimmersion, Fredholm, semi-Fredholm and regular point. Therefore the theorem presents several families $\mathcal {F}$ with $\dim M(x)\leq \infty $ and the trivial co-final set of $\mathcal {F}$ at $x_0$, all of which are $c^1$ integrable; while provides several kinds of solvable differential equations with an initial value. Let $B(E,F)$ be the set of all linear bounded
operators from $E$ into $F,\Phi_{m,n}$ the set of all double
splitting operators $T$ in $B(E,F)$ with $\dim N(T)=m$ and
$\mathrm{codim} R(T)=n$, and $F_k=\{T\in B(E,F):$ rank $T=k<\infty\}$.  Let
$\Lambda=B(E,F)\setminus\{0\}$ and $M(X)=\{T\in B(E,F):TN(X)\subset
R(X)\}$. In 1989, V. Cafagna in [Caf] introduced a geometrical method
for some partial differential equations and the family of subspaces
$\mathcal{F}=\{M(X):X\in\Lambda\}$. In this case the co-final set of $\mathcal {F}$ at $A$ is nontrivial in general (see the example in Section 4 ).                                     we prove the following theorems :  if $A\in\Lambda$ is double splitting, then
$\mathcal{F}$ at $A$ is smooth integrable; specially, for any one of
$F_k,\Phi_{m,n}(m,n<\infty),\Phi_{m,\infty}(m<\infty)$ and
$\Phi_{\infty,n}(n<\infty)$, write as $\Phi$, we have the global
result as follows, $\Phi$ is a smooth submanifold in $B(E,F)$ and
tangent to $M(X)$ at any $X\in\Phi$.

These results above seem to be useful for developing of differential
topology, global analysis, and the geometrical method in
differential equations (refer [Abr],[An] and [Caf]).
Let $U$ is an open set in $E,$ $S$ a $c^1$ submanifold in $U$, and $f$ a nonlinear functional from $U$ into $(-\infty ,\infty )$. A principle for critical point of $f$ under the constraint $S$ is given as follows, if $x_0 \in U$ is a critical point of $f$ under the constraint $S$,then $N( f'(x_0))\supset T_{x_0}S$ (see Theorem 4.3). The principle expands the principle in [Ma8] from the generalized regular constraint to the constraint of general $c^1$ submanifold in $E.$
 In view of the above Banach submanifolds with the explicit expression of the tangent space, the principle should be potential in application.\\

 Applying of generalized inverse is  usually neglected in pure
mathematical research. However, it plays an crucial role in this paper.
 We need the following theorems and concepts in generalized
inverse analysis  in Banach spaces.
Recall that $A^+\in B(F, E)$ is said to be a generalized inverse of
$A\in B(E,F)$\ provided $A^+AA^+=A^+$ and $A=AA^+A$; $A\in B(E,F)$
is said to be double splitting if $R(A)$ is closed and there exist
closed subspaces $R^+$ in $E$ and $N^+$ in $F$ such that
$E=N(A)\oplus R^+$ and $F=R(A)\oplus N^+$, respectively. It is well
known that $A$ has a generalized inverse $A^+\in B(F,E)$ if and only
if $A$ is double splitting (refer [Ma.1] and [N-C]).

Let $A\in B(E,F)$ be double splitting,
$A\not=0,$ and $A^+$ be a generalized inverse of $A$. Write
$V(A,A^+)=\{T\in
B(E,F):\|T-A\|<\|A^+\|^{-1}\},C_A(A^+,T)=I_F+(T-A)A^+,$ and
$D_A(A^+,T)=I_E+A^+(T-A).$ Then we have

{\bf Theorem 1.1}\quad {\it The following conditions for $T\in
V(A,A^+)$ are equivalent each other}:

(i) $R(T)\cap N(A^+)=\{0\}$;

(ii) $B=A^+C^{-1}_A(A^+,T)=D^{-1}_A(A^+,T)A^+$ {\it is the
generalized inverse of $T$ with $R(B)=R(A^+)$ and $N(B)=N(A^+)$};

(iii) $R(T)\oplus N(A^+)=F;$

(iv) $N(T)\oplus R(A^+)=E$;

(v) $(I_E-A^+A)N(T)=N(A)$;

(vi) $C^{-1}_A(A^+,T)TN(A)\subset R(A)$;

(vii) $R(C^{-1}_A(A^+,T)T)\subset R(A)$.

(Refer [N-C], [Ma1],[Ma.4] and [H-M].)

 Let $F_k=\{T\in B(E,F):$ rank$ T=k<\infty\}$, and  $\Phi_{m,n}=\{T\in B(E,F)$: $\dim N(T)=m<\infty$ and
 $\mathrm{codim} R(T)=n<\infty\}$.
 Let  $\Phi_{m,\infty}$ be the set of all semi-Fredholm operators
 $T$ with $\dim N(T)=m$ and $\mathrm{codim} R(T)=\infty$, and $\Phi_{\infty,n}$
 be the set of all semi-Fredholm operators $T$ with $\dim N(T)=\infty$
 and $\mathrm{codim} R(T)=n$. We have

 {\bf Theorem 1.2}\quad {\it Assume that  $A$  belongs to  any one of $F_k,
 \Phi_{m,n},\Phi_{m,\infty}$ and $\Phi_{\infty,n}.$ Then  the condition
 $R(T)\cap N(A^+)=\{0\}$ for $T\in V(A,A^+)$ holds if and only if
 $T$ and $A$ belong to the  same class.}
(Refer [N-C],  [Ma2])  and [Ma4].)

 {\bf Definition 1.1}\quad {\it
 Suppose that the  operator valued map
 $T_x$ from a topological space $X$ into $B(E,F)$ is continuous at
 $x_0\in X$, and  that $T_0=T_{x_0}$ is double splitting.
 $x_0$ is said to be a locally fine point of $T_x$
 provided there exist a generalized inverse $T^+_0$ of
 $T_0$ and a neighborhood $U_0$ (dependent on  $T^+_0)$ at $x_0$,
 such that
 $$R(T_x)\cap N(T^+_0)=\{0\},\quad\quad\forall x\in U_0.$$
By the way we state the reason why we define the locally fine point with the condition (i) in Theorem $1.1$
is its simple form and equivalence with the condition (ii) in Theorem $1.1.$ Essentially the show of the condition (ii) is a core property in generalized analysis, which is presented by Nashed, M.Z. and Chen ,X in [N-c]. We will refer to this property frequently.

{\bf Theorem 1.3}\quad{\it  The definition of  locally fine point
 $x_0$ of $T_x$ is independent of the choice of the generalized inverse
 $T^+_0$ of $T_0.$

{\bf Theorem 1.4}\quad (Operator rank theorem)\quad {\it Suppose
that the operator valued map $T_x:X\rightarrow B(E,F)$ is continuous
at $x_0\in X$ and $T_0$ is double splitting. Then the following
conclusion holds for arbitrary  generalized inverse $T^{+}_0$ of
$T_0$: there exists a neighborhood $U_0$ at $x_0$ such that $T_x$
has a generalized inverse $T^+_x$ for $x\in U_0$, and
$\lim\limits_{x\rightarrow x_0}T^+_x=T_0$, if and only if $x_0$ is a
locally fine point of $T_x$.}

In $1985,$ R. Penrose established the matrix rank theorem. Theorem $1.4$
expands the theorem from the case of matrices to  that of operators
in $ B(E,F)$ (refer [P] and [Ma $2$]).

{\bf Theorem 1.5}\quad {\it Suppose that $E_0$ and $E_1$ are two
closed subspaces in a Banach space $E$ with a common complement
$E_*$. Then there exists a unique operator $\alpha\in B(E_0,E_*)$
such that
$$E_1=\{e+\alpha e:\forall e\in E_0\}.$$
Conversely, $E_1=\{e+\alpha e:\forall e\in E_0\}$ for any $\alpha\in
B(E_0,E_*)$ is a closed subspace satisfying $E_1\oplus E_*=E$.}

For detail proofs of Theorems $1.1-1.5,$ see also the appendix in this
paper.

 \vskip 0.2cm\begin{center}{\bf 2\quad The Co-final Set and Frobenius
 Theorem}\end{center}\vskip 0.2cm

Let $E$ be a Banach space, and $\Lambda$ an open set in $E$. Assign
a subspace $M(x)$ in $E$ for every point $x$ in $\Lambda$,
especially, where the dimension of $M(x)$ may be infinite. In this
section, we consider the family $\cal{F}$ consisting of all $M(x)$
over $\Lambda,$ and  investigate the sufficient and necessary
condition for $\cal{F}$ being $c^1$ integrable at  a point $x$ in
$\Lambda$. For this purpose the concept of co-final set of $\mathcal {F}$ at $x_0 \in \Lambda$ is introduced.

{\bf Definition 2.1}\quad{\it Suppose $E=M(x_0)\oplus E_*$ for
$x_0\in\Lambda$. The set
$$J(x_0,E_*)=\{x\in\Lambda:M(x)\oplus E_*=E\},$$
is called a  co-final set of $\cal{F}$ at $x_0$.}

The next theorem tell us the co-final set and the integral
submanifold of $\cal{F}$ at $x_0$ are closely related.

{\bf Theorem 2.1}\quad{\it  If $\cal{F}$ is $c^1$ integrable at
$x_0\in\Lambda$, \mbox{and}  $S\subset E$ is the integral submanifold of
$\cal{F}$ at $x_0$, then there exist a closed subspace $E_*$ and a
neighborhood $U_0$ at $x_0$, such that
$$M(x)\oplus E_*=E,\quad\quad \forall x\in S\cap U_0,$$
i.e., $$J(x_0,E_*)\supset S\cap U_0.$$

{\bf Proof}\quad Recall that the  submanifold $S$ in the Banach
space $E$ is said to be tangent to $M(x)$ at $x\in S$ provided
$$M(x)=\{\dot{c}(0):\forall c^1\mbox{-curve }c(t)\subset S\mbox{ with }c(0)=x\}.
\eqno(2.1)$$
 By the definition of the submanifold in
Banach space, there exist a subspace $E_0$ splitting in $E$, say
$E=E_0\oplus E_1$, a neighborhood $U_0$ at $x_0$, and a $c^1$
diffeomorphism $\varphi:U_0\rightarrow\varphi(U_0)$ such that
$\varphi(S\cap U_0)$ is an open set in $E_0$. We claim
$$\varphi'(x)M(x)=E_0,\quad\quad\forall x\in S\cap U_0.\eqno(2.2)$$

Let $c(t)$ with $c(0)=x$ be an arbitrary $c^1$-curve contained in
$S\cap U_0$, then $\varphi'(x)\dot{c}(0)\in E_0$, and so
$\varphi'(x)M(x) \subset E_0$. Conversely, let
$r(t)=\varphi(x)+te\subset\varphi(S\cap U_0)$ for any $e\in E_0$,
and set $c(t)=\varphi^{-1}(r(t))\subset S\cap U_0$, then
$c(0)=\varphi^{-1}(\varphi(x))= x$ ,
and $\dot{c}(0)=(\varphi^{-1})'(\varphi(x)) e$
= $\varphi'(x)^{-1} e,$
so that $\varphi'(x)  M(x) \supset E_0.$ This proves that $(2.2)$
holds.

 Let  $E_x=\varphi'(x)^{-1}E_1$, and $E_*=E_{x_0}=\varphi '(x_0)^{-1} E_1.$ By
$(2.2) $ $ E= \varphi'(x_0)^{-1}(E_0\oplus E_1)=M(x_0) \oplus E_*. $
  In order to prove
 $$M(x)\oplus E_*=E\quad \quad \forall x\in S\cap U_0,$$
we consider the  projection  as follows,
 $$P_x=\varphi'(x)^{-1} \varphi '(x_0)P^{E_*}_{M(x_0)}\varphi'(x_0)^{-1}\varphi' (x) \quad\quad\forall
x\in S\cap U_0.\eqno(2.3)$$ Obviously, $P^2_x=P_x$, i.e., $P_x$ is a
projection on $E$. Next go to show
$$R(P_x)=M(x)\quad \mbox{and}\quad N(P_x)=E_x.$$
Indeed,
\begin{eqnarray*}  e \in N(P_x )& \Leftrightarrow &\varphi'(x_0)^{-1} \varphi'(x) e\in E_* \Leftrightarrow\varphi'(x)e\in \varphi'(x_0)E_* \\
&\Leftrightarrow & e\in\varphi'(x)^{-1}E_1 \, \, \mbox{because of} \, \, E_*=\varphi'(x_0)^{-1}E_1\\ &\Leftrightarrow &  e\in E_x,\end{eqnarray*}
and
 by $(2.3)$ and $(2.2),$
 $$R(P_x)= \varphi'(x)^{-1} \varphi' (x_0) M(x_0)= \varphi ' (x)^{-1}E_0=M(x) \, \, \mbox{for} \, \,
 x\in S\cap U_0 .$$ So $P_x=P^{E_x}_{M(x)}.$ Obviously $P_x$ is a generalized inverse of itself, and

 $$\lim\limits_{x\rightarrow
 x_0}P_x=P^{E_*}_{M(x_0)}.$$

 Let $X=S\cap U_0, T_x=P_x$, and  $T^+ _{x_0}=P^{E^*}_{M(x_0)}$. Then by Theorem $1.4$ there exists a neighborhood $V_0$ at $x_0$ in $S\cap U_0$, such that $\| P_x-P_{x_0}\| < \| P^{E_*}_{M(x_0)}\|$ and $R(P_x) \cap N(P_{x_0})= \{0\}$ for all $x\in V_0.$ Thus by the equivalence of the conditions (i) and (iii) in Theorem $1.1,$ $$ R(P_x) \oplus N( P_{M(x_0)}^{E_*})=E,\quad i.e.,\quad M(x)\oplus E_*\quad \forall x \in V_0.$$  For simplicity, still write $V_0$ as $U_0\cap S$. This shows $J(x_0,E_*)\supset S\cap U_0$ $\quad \Box$

{\bf Theorem 2.2}\quad{\it If $\mathcal{F}$ at $x_0$ is $c^1$
integrable, say that $S$ at $x_0$ is its integral submanifold in
$E$, then there exist a neighborhood $U_0$ at $x_0$,  a $c^1$
diffeomorphism $\varphi$ from $U_0$ onto  $\varphi(U_0)$ with
$\varphi'(x_0)=I_E$, and a subspace $E_*$ with $M(x_0)\oplus E_*=E,$ such that $V_0=\varphi(S\cap U_0)$ is an open
set in $M(x_0)$,
$$J(x_0,E_*)\supset S\cap U_0 \quad  \mbox{and} \quad  \varphi'(x) M(x)=M(x_0),\quad\quad\forall x\in {S\cap U_0}.\eqno(2.4)$$

{\bf Proof} The following conclusion has been indicated in the proof of Theorem $2.1 :$  there are the subspace $E_0$ with $E=E_0\oplus E_1,$ the neighborhood $U_0$ at $x_0,$  and the
diffeomorphism $\varphi $ from $U_0$ onto $\varphi (U_0)$ such that $ \varphi (S\cap U_0)$ is an open set $E_0$, $\varphi '(x)M(x)=E_0 \, \mbox {and} \, M(x)\oplus E_*=E \,  \forall x\in {S\cap U_0},$ where $ E_*= \varphi '(x_0)^{-1}E_1.$
Consider
$$\varphi_1(x)=\varphi'(x_0)^{-1}\varphi(x),\quad\quad\forall x\in
S\cap U_0.$$ Obviously $\varphi'_1(x_0)=I_E$.
 Note $\varphi'(x_0)^{-1}E_0=M(x_0)$.
\begin{eqnarray*}
\varphi'_1(x)M(x)&=&\varphi'(x_0)^{-1}\varphi'(x)M(x)\\
&=&\varphi'(x_0)^{-1}E_0=M(x_0) \quad \quad \forall x\in S\cap U_0,
\end{eqnarray*}
and
$$\varphi_1(S\cap U_0)=\varphi'(x_0)^{-1}(\varphi(S\cap U_0))$$
is an open set in $M(x_0)$ because of $\varphi(S\cap U_0)$ being open an
set in $E_0$. Finally write $\varphi_1$ as $\varphi$.
 The theorem is proved. \quad $\Box$

 According to Theorem $1.5,$  $M(x)$ for $x\in
J(x_0,E_*)$ has the coordinate expression
$$M(x)=\{e+\alpha (x)e:\forall e\in M (x_0)\}$$
where $\alpha\in B(M(x_0),E_*).$ This shows that $M(x)$ for $x \in J(x_0,E_*)$ possesses an unique operator valued coordinate $\alpha (x)\in B(M(x_0),E_*).$

We now state the  Frobenius theorem in Banach spaces. For simplicity, write $M_0=M(x_0)$ in the sequel.

{\bf Theorem 2.3}\ (Frobenius theorem)\quad{\it $\cal{F}$ is $c^1$
integrable at $x_0$ if and only if the following conditions hold:}

(i) $M_0$  splits in $E$, say $E=M_0\oplus E_*$;

(ii) there exist a neighborhood} $V$ at $ P^{E_*}_{M_0}x_0$
in $M_0$, and a $c^1$  map $\psi:V\rightarrow E_*$,
such that $x+\psi(x)\in J(x_0,E_*)$  for all $x\in V$,
and $\alpha(x+\psi(x))$  is continuous in $V$;

(iii) $\psi$  satisfies
$$\begin{array}{rl}
&\psi'(x)=\alpha(x+\psi(x))\quad \quad \mbox{ for\ all} \quad x\in V,\\
&\psi(P^{E_*}_{M_0}x_0)=P^{M_0}_{E_*}x_0,\end{array}\eqno(2.4)$$ where  $\psi'(x)$ is the Fr\'{e}chet derivative of  $\psi$
$at$ $x$.

{\bf Proof}\quad Assume that $\cal{F}$ is $c^1$-integrable at $x_0$.
Go to prove that the conditions (i), (ii) and (iii) in the theorem
hold. The condition (i) follows from Theorem $ 2.1$. Next go to prove
that the conditions (ii) and (iii) in the theorem hold.
Consider the map $\Gamma$ as follows,
$$\Gamma(x)=\varphi(x)+P^{E_*}_{M_0}x_0-\varphi(x_0)\quad\quad\forall
x\in U_0,$$ where $\varphi $ and $U_0$  are as indicated in Theorem $2.2$. It is obvious that $\Gamma$ from $U_0$ onto the open set
 $\Gamma(U_0)=\varphi(U_0)+P^{E_*}_{M_0}x_0-\varphi(x_0)$ is $c^1$
diffeomorphism. Write the open set $\varphi(S\cap
U_0)+P^{E_*}_{M_0}x_0-\varphi(x_0)$ in $M_0$ as $V_1$. Then by Theorem $2.2$
$$\Gamma^{-1}(V_1)=\Gamma^{-1}(\Gamma(S\cap U_0))=S\cap U_0\subset
J(x_0,E_*).$$ Let $\varphi_0=P^{E_*}_{M_0}\Gamma^{-1}  \mbox{and}
 \, \varphi_1=P^{M_0}_{E_*}\Gamma^{-1}. $ Note $\Gamma(x_0)=P^{E_*}_{M_0}x_0  \, \, \mbox{and} \, \, \Gamma'(x_0)=\varphi'(x_0)=I_E$.
Directly
\begin{eqnarray*}
\varphi'_0(P^{E_*}_{M_0}x_0)&=&P^{E_*}_{M_0}(\Gamma^{-1})'(P^{E_*}_{M_0}x_0)\\
&=&P^{E_*}_{M_0}(\Gamma^{-1})'(\Gamma(x_0))\\
&=&P^{E_*}_{M_0}\Gamma'(x_0)^{-1}\\
&=&P^{E_*}_{M_0}\varphi'(x_0)^{-1}=P^{E_*}_{M_0}
\end{eqnarray*}
\mbox{and}
$$\varphi_0(P^{E_*}_{M_0}x_0)=P^{E_*}_{M_0}\Gamma^{-1}(P^{E_*}_{M_0}x_0)=P^{E_*}_{M_0}x_0.$$
By the inverse mapping theorem  there exists a neighborhood $V_*$ at $P^{E_*}_{M_0}x_0$ in
$V_1\subset M_0$ such that $\varphi_0$ from $V_*$ onto
$V=\varphi_0(V_*)$ is $c^1$ diffeomorphism, where $V$ is an open set
in $ M_0 $   \mbox{and}  contains $P^{E_*}_{M_0} x_0$ since $\varphi_0 (V_*)\subset M_0 \, \, \mbox{and}  \, \, P_{M_0}^{E_*}$ is a fixed point of $\varphi _0.$
Let $\psi=\varphi_1\circ\varphi^{-1}_0$ and $y=\varphi _0(x) $ for $  x \in V_*.$ Since $\Gamma ^{-1} (V_1 ) = S\cap U_0  \in J(x_0,E_*)$
$$\Gamma^{-1}(x)=\varphi_0(x)+\varphi_1(x)=y+\psi(y) \in J(x_0,E_*) \, \mbox{ for \ all } \, y \in V .$$ ( It is good in the sequel to bear $V=\varphi_0(V_*)\subset M_0$ in mind.)
According to Theorem $1.5$ we have the operator $\alpha\in B(M_0,E_*)$
satisfying
$$M(\Gamma^{-1}(x))=\{e+\alpha(y+\psi (y))e:\forall e\in M_0\}.$$
 Since $\Gamma'
(\Gamma^{-1}(x))=\varphi'(\Gamma^{-1}(x))$ and Theorem $2.2$
\begin{eqnarray*}
M(\Gamma^{-1}(x))&=&\varphi'(\Gamma^{-1}(x))M_0\\
&=&\Gamma'(\Gamma^{-1}(x))^{-1}M_0)\\
&=&(\Gamma^{-1})'(x)M_0\\
&=&(\varphi'_0(x)+\varphi'_1(x))M_0.
\end{eqnarray*}
Then, for any $e\in M_0$ there is $e_*\in M_0$ such that
$e+\alpha(y+\psi(y))e=\varphi'_0 ( x )e_*+\varphi'_1(x)e_*,$
and so,
$e=\varphi'_0(x)e_* \, \mbox{and} \, \ \alpha(y+\psi(y))e=\varphi'_1(x)e_*.$ Thus
\begin{eqnarray*}
\alpha(y+\psi(y))e&=&\varphi'_1(x)\varphi'_0(x)^{-1}e\\
&=&\varphi'_1(x)(\varphi^{-1}_0)'(y)e\\
&=&(\varphi_1\circ\varphi^{-1}_0)'(y)e\\
&=&\psi'(y)e, \end{eqnarray*} and so
$\psi'(y)=\alpha(y+\psi(y)) \, \forall y\in V.$ Hereby the
conditions (ii) and  (iii) follow from $\psi$ is of $c^1$. Note
$\varphi_0(P^{E_*}_{M_0}x_0)=P^{E_*}_{M_0}x_0 \,  \mbox{and}  \,
\Gamma(x_0)=P^{E_*}_{M_0}x_0.$ Clearly
\begin{eqnarray*}
\psi(P^{E_*}_{M_0}x_0)&=&\varphi_1(P^{E_*}_{M_0}x_0)\\
&=&P^{M_0}_{E_*}\Gamma^{-1}(P^{E_*}_{M_0}x_0)\\
&=&P^{M_0}_{E_*}x_0.
\end{eqnarray*}
Now, the necessity of the theorem is proved.

Assume that the conditions (i), (ii) and (iii) hold. Go to show that
$\cal{F}$ is $c^1$ integrable at $x_0$. Let $S=\{x+\psi(x):\forall
x\in V\}$,
$$V^*=\{x\in E:P^{E_*}_{M_0}x\in V\}(\supset V) \quad \mbox{ and} \quad
\Phi(x)=x+\psi(P^{E_*}_{M_0}x) \,  \, \forall x\in V^*.$$
 Where $V, M_0,\psi$ and $E_*$ are as indicated in the conditions
 (i), (ii) and (iii).
Obviously, $V^*$ is an open set in $E$, $\Phi(V)=S,$ and $\Phi$ is a
$c^1$ map. Moreover, we are going to prove that
$\Phi:V^*\rightarrow\Phi(V^*)$ is a diffeomorphism. Evidently, if
$\Phi(x_1)=\Phi(x_2)$ for $x_1,x_2\in V^*$,  then
$$P^{E_*}_{M_0}(x_1-x_2)+\psi(P^{E_*}_{M_0}x_1)-\psi(P^{E_*}_{M_0}x_2)+P^{M_0}_{E_*}(x_1-x_2)=0$$
and so
$$P^{E_*}_{M_0}x_1=P^{E_*}_{M_0}x_2 \ \mbox{ and \ hence} \
P^{M_0}_{E_*}x_1=P^{M_0}_{E_*}x_2.$$
 This shows that
$\Phi:V^*\rightarrow \Phi(V^*)$ is one-to-one. Now, in order to show
that $\Phi:V^*\rightarrow\Phi(V^*)$ is a $c^1$ diffeomorphism, we
need only to show that $\Phi(V^*)$ is an open set in $E$. By the
inverse map theorem it is enough to examine that $\Phi'(x)$ for any
$x\in V^*$ is invertible in $B(E)$. According to the condition (iii)
in the theorem
$$\begin{array}{rllr}
\Phi'(x)&=&I_E+\psi'(P^{E_*}_{M_0}x)P^{E_*}_{M_0}\\
&=&P^{E_*}_{M_0}+\alpha(P^{E_*}_{M_0}x+\psi(P^{E_*}_{M_0}x))P^{E_*}_{M_0}+P^{M_0}_{E_*}\end{array}\eqno(2.5)$$
for any $x\in V^*$. We claim $N(\Phi'(x))=\{0\}$. If $\Phi'(x)e=0$,
i.e.,
$P^{E_*}_{M_0}e+\psi'(P^{E_*}_{M_0}x)P^{E_*}_{M_0}e+P^{M_0}_{E_*}e=0,$
then $P^{E_*}_{M_0}e=0$ and so, $P^{M_0}_{E_*}e=0$. This says
$N(\Phi'(x))=\{0\}$ for any $x\in V^*$. Next go to  verify that
$\Phi'(x)$ is surjective.
For abbreviation, write $M(y)=M_*$ and $y=P^{E_*}_{M_0}x+\psi(P^{E_*}_{M_0}x)$ for any $x\in V^*$.
Obviously, $y \in S$ and by the assumption (ii), $y \in J(x_0,E_*)$,
i.e., $M_*\oplus E_*=E$. Hence by Theorem $1.5$ we have  for any $e\in
E$ there exist $e_0\in M_0$ and $\alpha\in B(M_0,E_*)$ such that
$P^{E_*}_{M_*}e=e_0+\alpha(y)e_0$. Set $e_*=e_0+P^{M_*}_{E_*}e$,
then by $(2.5),$
$\Phi'(x)e_*=e_0+\alpha(y)e_0+P^{M_*}_{E_*}e=P^{E_*}_{M_*}e+P^{M_*}_{E_*}e=e.$
 This says that $\Phi'(x)$ for any $x\in V^*$ is surjective. Now it is
 proved that $\Phi^{-1}$  is a
 $c^1$-diffeomorphism  from the open set $\Phi(V^*)$ onto $V^*$, and that $\Phi^{-1}(S)=V$ is an open set in $M_0$.
  Then by the condition (i), $S$ is a $c^1$ submanifold in $E$. Finally go to show that
  $S$ is tangent to $M(x)$ at any point $x\in S$.
 Write
  \[T(x+\psi(x))=\{\dot{c}(0):\forall c^1\mbox{-curve }c(t)\subset
  S \, \mbox{with}\ c(0)=x+\psi(x)\}\]  for any $x\in V$. Namely we have to examine  $  M(y)=T(y)$ for $y  \in S.$
By $S=\{y: y=x+\psi(x) \, \forall x\in V\}$, $P_{M_0}^{E_*}y=x$ and so, $y=P_{M_0}^{E_*}y + \psi (P_{M_0}^{E_*}y)$ for any $y\in S.$ Hence we have for any $c^1$-curve $\gamma (t)\subset S $ with $\gamma (0)=y,$  $$\gamma (t)=P_{M_0}^{E_*} \gamma (t) + \psi (P_{M_0}^{E_*} \gamma (t)).$$
Clearly $P_{M_0}^{E_*} \gamma(t)$ is a $c^1$-curve $\subset V$ with $P_{M_0}^{E_*} \gamma(0)=P_{M_0}^{E_*}y,$ write it as $c(t)$. By the condition(iii) \begin{eqnarray*}
\dot{\gamma}(0)&=&\dot {c}(0)+\psi '(P_{M_0}^{E_*}y)\dot{c}(0)\\&=& \dot{c}(0)+\alpha (P_{M_0}^{E_*}y+\psi (P_{M_0}^{E_*}y))\dot{c}(0)\\&=&\dot{c}(0)+\alpha (y)\dot{c}(0)+\alpha(y)\dot{c}(0)\in M(y) \end{eqnarray*}because of $\dot{c}(0)\in M_0.$ This shows $M(y)\supset T(y)$ for any $y \in S.$ Conversely assume $e\in M(y),$ say $e=e_0+\alpha (y)e_0$ for some $e_0\in M_0.$ Consider the $c^1$-curve $c(t)$ as follows, $c(t)=P_{M_0}^{E_*}y+te_0 $ for $t$ near 0. Then $\gamma(t)=c(t)+\psi (c(t))\subset S$ and $\gamma(0)=y;$ while by the condition (iii), \begin{eqnarray*} \dot{\gamma}(0)&=& e_0+ \psi'(P_{M_0}^{E_*}y)e_0 \\ &=& e_0+\alpha (P_{M_0}^{E_*}y +\psi (P_{M_0}^{E_*}y))e_0\\ &=& e_0+\alpha (y))=e\end {eqnarray*} for $y \in S.$ This shows $T(y)\supset M(y).$ Therefore $M(y)=T(y)$ for all $y \in S.$ The proof ends.\quad$\Box$

In the sequel, it will be seen that the co-final set and the
coordinate operator $\alpha(x)$ for $M(x)$ are essential to
Frobenious  theorem.

\vskip 0.2cm\begin{center}{\bf 3\quad $C^1$ Integrable Family of
Subspaces with Trivial Co-final Set }\end{center}\vskip 0.2cm

In this section we will discuss $c^1$ integrable family of subspaces
with trivial co-final set.

\textbf{Definition $3.1$}\quad {\it The co-final set $J(x_0,E_*)$ is said
to be trivial provided that $x_0$ is an inner point of $J(x_0,E_*)$,
and $\alpha(x)$ is continuous in some neighborhood at  $x_0$.

\textbf { Theorem $3.1$} \quad {\it If $J(x_0,E_*)$ is trivial, then the $c^1$ integrability of $\mathcal{F}$ at $x_0$ and the solvability of the initial problem $(2.4)$ are equivalent.}

$\quad ${\bf Proof}\quad By Definition $3.1$ there exists a neighborhood $W_0$ at $x_0$ such that $W_0\subset J(x_0,E_*).$ Without loss of generality, we assume that $\alpha $ is continuous in $W_0.$ In what follows, we claim that there  exist a neighborhood $V$ at $P_{E_*}^{M_0}x_0$ in $M_0$ and a neighborhood $U$ at $P_{E_*}^{M_0} x_0$ in $E_*$ such that the condition (ii) in Theorem $2.3$ for any $c^1$ map $\psi$ from $V$ into $U$ with $\psi (P^{E_*}_{M_0} x_0)=P_{E_*}^{M_0} x_0 $ holds. Consider the continuous map $x+y$ from $M_0\times E_*$ into $E_*.$ Clearly,  there exist a neighborhood  at $P^{E_*}_{M_0} x_0$ in $M_0$, still write it as $V$ , and a neighborhood  at $P_{E_*}^{M_0} x_0$ in $E_*$, still write it as $U$, such that $x+y\in W_0 $ for $ x \in V $ and  $ y\in U.$  Thus  $x+\psi(x)\in W_0 \  \forall x\in V$ and $\alpha(x+\psi (x))$ is continuous in $V$. This says  that  the conditions (ii) holds for arbitrary $c^1$ map $\psi$ with $\psi (P^{E_*}_{M_0}x_0)=P_{E_*}^{M_0} x_0 $. (This is why we call  $J(x_0,E_*)$ in Definition $3.1$ to be trivial.)  Finally by Theorem $2.3$  the theorem is proved.\quad $\Box$\\
The next example can be illustrate Theorem $3.1$ although it simple.

 {\bf Example}\quad{\it Let $E=\mathbf{R}^2,
  \Lambda=\mathbf{R}^2\setminus(0,0)$ and
  $$M(x,y)=\{(X,Y)\in\mathbf{R}^2:Xx+Yy=0\},\quad\quad\forall
  (x,y)\in\Lambda.$$
  Consider the family of subspaces $\mathcal{F}=\{M(x,y):\forall
  (x,y)\in \Lambda\}$. Apply Frobenius theorem in Banach spaces to
  determine the  integral curve of $\mathcal{F}$ at $(0,1)$.}

  Set $U_0=\{(x,y)\in\mathbf{R}^2:y>0\}$ and
  $E_*=\{(0,y)\in\mathbf{R}^2:\forall y\in R\}$.
  Obviously, $U_0\subset \Lambda$, and
  $$M(x,y)\oplus E_*=\mathbf{R}^2\quad\quad\forall (x,y)\in U_0$$
  since
  $M(x,y)\cap E_*= \{(0,0)\} \, \forall(x,y)\in U_0.$

  Hence $J((0,1),E_*)\supset U_0$. This shows that $(0,1)$ is an inner point of  $J((0,1),E_*)$.
  Next go to determine $\alpha$ in the equation $(2.4).$ Note
  $M_0=M(0,1)=\{(X,0):\forall X\in R\}$. Evidently,
  \begin{eqnarray*}
  M(x,y)&=&\{(X,-\frac{x}{y}X):\forall X\in R\}\\
  &=&\{(X,0)+(0,-\frac{x}{y}X):\forall X\in
  R\}\end{eqnarray*}
  for all $(x,y)\in U_0$.
   Theorem $1.5$ shows that the operator $\alpha \in B(M_0,E_*)$ is unique, so we conclude
  $$\alpha(x,y)(X,0)=(0,-\frac{x}{y}X)\quad\quad\forall (x,y)\in
  U_0.$$
  We claim that $\alpha:U_0\rightarrow B(M_0,E_*)$ is
  continuous in $U_0$. Obviously
  \begin{eqnarray*}
  &&\left\|(\alpha(x+\Delta x,y+\Delta
  y)-\alpha(x,y))(X,0)\right\|\\
  &&\quad\quad=\left\|(0, (\frac{x}{y}-\frac{x+\Delta
  x}{y+\Delta y})X)\right\|=\left|\frac{x+\Delta x}{y+\Delta
  y}-\frac{x}{y}\right|\|(X,0)\|
  \end{eqnarray*}
  for any $(x,y)\in U_0$, where $\|,\|$ denotes the norm in
  $\mathbf{R}^2$.

  So
  $$\|\alpha(x+\Delta x,y+\Delta
  y)-\alpha(x,y)\|=\left|\frac{x+\Delta x}{y+\Delta
  y}-\frac{x}{y}\right|,$$
  where $\|,\|$ denotes the norm in $B(M_0,E_*).$ Hereby, one can
  conclude that $\alpha:U_0\rightarrow B(M_0,E_*)$ is continuous in
  $U_0$. This infers that $J((0,1),E_*)$ is trivial. By Theorem $3.1$ we now ought to solve $(2.4).$ Let
  $$V=\{(x,0):|x|<1\}\subset M((0,1))$$
  and
  $$\psi((x,0))=(0,y(x))\in E_*\quad{\mbox for \ any} \
  (x,0)\in V,$$
  where $y(x):(-1, 1)\mapsto {\bf R}$ is a $c^1$ function. Then by the equation $(2.4).$
 \begin{eqnarray*}
&&\frac{dy}{dx}=-\frac{x}{y},\quad\quad\forall x\in (-1,1),\\
&&y(0)=1.\end{eqnarray*} It is easy to see that the solution is
$y=\sqrt{1-x^2}$ for all $x\in(-1,1)$. So, $S=\{(x, \sqrt{1-x^2}):
x\in(-1, 1)\}$ is the smooth integral curve of $\mathcal F=\{M(x,
y): (x, y)\in\Lambda\}$ at $(0,1)$. (Refer [Ma $4$] also [Ma$7$]).

Let $f(x)$ be a $c^1$ map from an open set $U\subset E$ into $F,
N_x=N(f'(x))$ for $x\in U$ \, \mbox{and} \, $\mathcal{F}=\{N_x:\forall x\in U\}.$
Recall that $x_0\in U$ is said to be a generalized regular point of
$f$ provided $x_0$ is a locally fine point of $f'(x)$ ( for details
see the appendix in this paper, and also [Ma$3$], [Ma$4$] and [Ma$8$]).
 When $x_0$ is a generalized regular point of $f$, by Definition $1.1$ and Theorem $1.3$ we have for any generalized inverse
$T^+_0$ of $f'(x_0)$ there exists a neighborhood $U_0$ at $x_0$ such
that
$$R(f'(x))\cap N(T^+_0)=\{0\},\quad\forall x\in U_0.$$
Since $f'(x)$ is continuous at $x_0$, one can assume
$$\|f'(x)-f'(x_0)\|<\|T^{+}_0\|^{-1}\quad \mbox{and}\quad R(f'(x))\cap N(T^+_0)=\{0\}\eqno(3.1)$$ for all $x\in U_0$. Let
$E_*=R(T^+_0)$.

\textbf{Theorem $3.2$}\quad{\it If $x_0$ is a generalized regular
point of $f$ then $J(x_0,E_*)$ is trivial, and
$\mathcal{F}=\{N(f'(x)):\forall x\in U_0\}$ is $c^1$ integrable at
$x_0$.}

\textbf{Proof}\quad We claim that $x_0$ is an inner point of
$J(x_0,E_*)$. By $(3.1),$
$$\|f'(x)-f'(x_0)\|<\|T^+_0\|^{-1},\quad \forall x\in U_0 ,$$
one can apply the equivalence of the conditions (i) and (iv) in
Theorem $ 1.1$ to $f'(x)$ for any $x\in U_0$, so that
$$N(f'(x))\oplus E_*=E\quad\forall x\in U_0,$$
i.e., $J(x_0,E_*)\supset U_0$ . So $x_0$ is a inner point of
$J(x_0,E_*)$. Next go to verify that $\alpha (x)$ is continuous in $U_0.$ For this we will
need to investigate the coordinate operator $\alpha(x)$ of $N_x$ more
carefully. The generalized point of $f$, $x_0$ bears such a explicit expression of $\alpha (x)$ for $x\in U_0$ with $f'(x)$  that it follows that $\alpha (x)$ is continuous in $U_0$ from the expression.

Evidently
$$P^{E_*}_{N_x}P^{E_*}_{N_0}e=P^{E_*}_{N_x}(P^{E_*}_{N_0}e+P^{N_0}_{E_*}e)=P^{E_*}_{N_x}e=e,\quad \forall
e\in N_x$$ and
$$P^{E_*}_{N_0}P^{E_*}_{N_x}x=P^{E_*}_{N_0}(P^{E_*}_{N_x}x+P^{N_x}_{E_*}x)=P^{E_*}_{N_0}x=x,\quad  \forall
x\in N_0,$$ where $N_0=N(f'(x_0)).$  So
$$e=P^{E_*}_{N_0}e+P^{N_0}_{E_*}e=P^{E_*}_{N_0}e+P^{N_0}_{E_*}P^{E_*}_{N_x}P^{E_*}_{N_0}e,$$
for all $e\in N_x$.

Since $\alpha(x)$ in Theorem $1.5$ is unique, we conclude
$$ \left.\alpha(x)=P^{N_0}_{E_*}P^{E_*}_{N_x}\right|_{N_0}\ \ \ \ \forall x\in U_0.$$

Take the places of $T$ and $A$ in Theorem $1.1$ by $T_x=f'(x)$ for any
$x\in U_0$  and $T_0=f'(x_0)$, respectively.  Due to $(3.1)$ one can apply Theorem $1.1$ to $f'(x)$ for $x\in U_0$. Then  the
equivalence of the conditions (i) and (ii) in Theorem $1.1$ shows that
$T_x$ has the generalized inverse $T^+_x$ as follows,
$ T^+_x=D^{-1}_{T_0}(T^+_0,T_x)T^+_0 $ with $ N(T^+_x)=N(T^+_0) $ and $ R(T^+_x)=R(T^+_0) \  \forall x \in U_0.$
Thus, by the definition of generalized inverse we have
$T^+_xT_x=P^{N_x}_{R(T^+_x)}$  and
\begin{eqnarray*}
 P^{E_*}_{N_x}&(=&P^{R(T^+_0)}_{N_x})=P^{R(T^+_x)}_{N_x}
 =I_E-T^+_xT_x \\ &=&I_E-D^{-1}_{T_0}(T^+_0,T_x)T^+_0T_x.
\end{eqnarray*}
So
\begin{eqnarray*}
\alpha(x)&=&\left.P^{N_0}_{E_*}P^{E_*}_{N_x}\right|_{N_0}\\
&=&P^{N_0}_{E_*}(I_E-D^{-1}_{T_0}(T^+_0,T_x)T^+_0T_x)|_{N_0}\\
&=&P^{N_0}_{E_*}\left.\left\{I_E-D^{-1}_{T_0}(T^+_0,T_x)(P^{E_*}_{N_0}+T^+_0T_x-P^{E_*}_{N_0})\right\}\right|_{N_0}\\
&=&P^{N_0}_{E_*}D^{-1}_{T_0}(T^+_0,T_x)P^{E_*}_{N_0}
\end{eqnarray*}
for all $x\in U_0.$ Then it follows that $J(x_0,E_*)$ is trivial.
Now to end the proof, by Theorem $3.1$ it is enough to verify that the problem $(2.4)$ is solvable. Consider
$$\varphi(x)=T^+_0(f(x)-f(x_0))+P^{E_*}_{N_0}x \ \ \  \forall x\in U_0. $$
Obviously
$\varphi'(x)=P^{E_*}_{N_0}+T^+_0T_x=D_{T_0}(T^+_0,T_x),\quad
\varphi(x_0)=P^{E_*}_{N_0}x_0$
and $\varphi'(x_0)=I_E$.
So there exists a neighborhood at $x_0$, without loss of generality,
still write it as $U_0$, such that $\varphi$ from $U_0$ onto
$\varphi(U_0)$ is $c^1$ diffeomorphism.

Let $y=\varphi(x)$ for $x\in U_0$. Note
$(\varphi'(x))^{-1}=(\varphi^{-1})'(y).$
Then
\begin{eqnarray*}
\alpha(x)&=&P^{N_0}_{E_*}D^{-1}_{T_0}(T^+_0,T_x)|_{N_0}\\
&=&P^{E_*}_{N_0}(\varphi'(x))^{-1}|_{N_0}\\
&=&P^{E_*}_{N_0}(\varphi^{-1})'(y)|_{N_0} \end{eqnarray*}
and so,
$$\alpha(\varphi^{-1}(y))=P^{E_*}_{N_0}(\varphi^{-1})'(y)|_{N_0}\eqno(3.2)$$
for $y\in \varphi(U_0)$.
Comparing the equation $(2.4)$ with the equality $(3.2)$ one can comprehend an approach with $\varphi^{-1}$ to get the solution $\psi$ of the problem $(2.4).$ Next go to realize this ideal.
Let $V_*=\varphi(U_0)\cap N_0$,
$$\varphi_0(y)=P^{E_*}_{N_0}\varphi^{-1}(y)\mbox{ and }
\varphi_1(y)=P^{N_0}_{E_*}\varphi^{-1}(y)\quad \forall y\in V_*.$$ We claim that $\varphi_0$ is a $c^1$
diffeomorphism from $V_*$ onto $\varphi_0(V_*)$.
Obviously
$\varphi_0 $  maps $ V_*(\subset N_0)$  into $ N_0$
and
\begin{eqnarray*}
\varphi'_0(P^{E_*}_{N_0}x_0)&=&P^{E_*}_{N_0}(\varphi^{-1})'(P^{E_*}_{N_0}x_0)\\
&=&P^{E_*}_{N_0}(\varphi'(x_0))^{-1}\\
&=&P^{E_*}_{N_0}.
\end{eqnarray*}
So, there exists a neighborhood at $P^{E_*}_{N_0}x_0$ in $N_0$,
without loss of generality still write it as $V_*$, such that
$\varphi_0$ from $V_*$ onto $\varphi_0(V_*)$ is $c^1$ diffeomorphism.
Moreover go to verify
$$\varphi'_0(y)=P^{E_*}_{N_0} \quad\forall y\in V_*.\eqno(3.3)$$ Note the following two equalities for $x\in U :$
$$
(\varphi'(x))^{-1}=D^{-1}_{T_0}(T^+_0,T_x) \quad\forall x\in U_0$$
and
$$P^{E_*}_{N_0}D^{-1}_{T_0}(T^+_0,T_x)=I_E-T^+_0T_xD^{-1}_{T_0}(T^+_0,T_x) \quad \forall
y\in V_* .$$

Then
\begin{eqnarray*}
\varphi'_0(y)&=&P^{E_*}_{N_0}(\varphi^{-1})'(y)\\
&=&P^{E_*}_{N_0}(\varphi'(x))^{-1}\\
&=&P^{E_*}_{N_0}D^{-1}_{T_0}(T^+_0,T_x)\\
&=&P^{E_*}_{N_0}P^{E_*}_{N_0}D^{-1}_{T_0}(T^+_0,T_x)\\
&=&P^{E_*}_{N_0}(I_E-T^+_0T_x D^{-1}_{T_0}( T^+_0,T_x))\\
&=&P^{E_*}_{N_0}\end{eqnarray*} because of
$R(T_0^+)=E_*,$ where $x=\varphi_0(y)\in N_0.$

Let $V=\varphi_0(V_*)$. Clearly, $P^{E_*}_{N_0}x_0\in V$ and $V$ is an
open set in $N_0$.

Finally go to prove that
$\psi(z)=(\varphi_1\circ\varphi^{-1}_0)(z)$ for all $z\in V$ is the
solution of $(2.4).$

By $(3.3)$
\begin{eqnarray*}
\psi'(z)|_{N_0}&=&\varphi'_1(\varphi^{-1}_0(z))(\varphi^{-1}_0)'(z)|_{N_0}\\
&=&\varphi'_1(y)(\varphi'_0(y))^{-1}|_{N_0}\\
&=&P^{N_0}_{E_*}(\varphi^{-1})'(y)P^{E_*}_{N_0}|_{N_0}\end{eqnarray*}
where $y=\varphi^{-1}_0(z)\in V_*.$

Since $\varphi^{-1}_0(V)=V_*\subset\varphi(U_0),$
$y=\varphi^{-1}_0(z)\in\varphi(U_0)$ for all $z\in V.$
Then by $(3.2)$
 $\psi'(z)|_{N_0}=\alpha(\varphi^{-1}(y))|_{N_0}.$
Meanwhile
 \begin{eqnarray*}
\varphi^{-1}(y)&=&\varphi_0(y)+\varphi_1(y)\\
&=&z+(\varphi_1\circ\varphi^{-1}_0)(z)\\
&=&z+\psi(z). \end{eqnarray*}
Therefore
$\alpha(z+\psi(z))=\psi'(z)$ for all $ z\in V.$
Note $\varphi (x_0)=P^{E_*}_{N_0}x_0 $ and $ \varphi_0(P^{E_*}_{N_0}x_0)=P^{E_*}_{N_0}x_0.$ Clearly
\begin{eqnarray*}
\psi(P^{E_*}_{N_0}x_0)&=&\varphi_1\left(\varphi^{-1}_0(P^{E_*}_{N_0}x_0)\right)\\
&=&\varphi_1(P^{E_*}_{N_0}x_0)=P_{E_*}^{N_0}\varphi^{-1}(P^{E_*}_{N_0}x_0)\\
&=&P^{N_0}_{E_*}x_0\end{eqnarray*}
By Theorem $3.1$  $\mathcal{F}$ at $x_0$ is
$c^1$ integrable.\quad$\Box$

\textbf{Remark}

( i ) The generalized regular point is often the case in non linear functional analysis as indicated in Section $1,$and so does the $c^1$ integrable family $\mathcal{F}$ at a generalized regular point $x_o$ with $dimM(x)\leq \infty$ and trivial co-final set $J(x_0,E_*)$;

(ii) Theorem $3.2$ provides several kinds of solvable differential equations with an initial value;

(iii) for more of applications of generalized regular point see 6 in the appendix of this paper.

\vskip 0.2cm\begin{center}{\bf 4\quad A Family of Subspaces with
Non-trivial Co-final Set and Its Smooth Integral
Submanifolds}\end{center}\vskip 0.2cm

Let  $\Lambda =B(E,F)\setminus\{0\}$ and $M(X)=\{T\in
B(E,F):TN(X)\subset R(X)\}$ for $X\in\Lambda$. V. Cafagna introduced
the  geometrical method for some partial differential equations and
presented the family of subspaces $\mathcal{F}=\{M(X)\}_{X\in\Lambda}$
in [Caf]. Now we take the following example to  illustrate that when
A is neither left nor right invertible in $B(R^2), $ the co-final
set of $\cal{F}$ at $A$ is non-trivial.

{\bf Example}\quad Consider the space $B(R^2)$  consisting of
all real $2\times 2$ matrices. Let $A$ and $\mathbf{E}_*$ be
$$\left(\begin{array}{cc}1&0\\ 0&0\end{array}\right)\quad{\rm
and}\quad\left\{\left(\begin{array}{cc}0&0\\
0&t\end{array}\right):\forall t\in R\right\},$$ respectively.
Obviously \begin{eqnarray*} &&N(A)=\{(0,x):\forall x\in R\},\ \
R(A)=\{(x,0),\forall x\in R\},\\
&&M(A)=\left\{\left(\begin{array}{cc}t_{11}&t_{12}\\
t_{21}&0\end{array}\right):\forall t_{11},t_{12}\ \mbox{ and}\
t_{21}\in R\right\},
\end{eqnarray*}
\noindent and\\[-10pt]
$$\hspace*{-130pt}M(A)\oplus\mathbf{E}_*=B(R^2).$$
To show that $J(A,\mathbf{E}_*)$ is nontrivial. Consider
$$A_\varepsilon=\left(\begin{array}{cc}1&0\\
0&\varepsilon\end{array}\right),\quad\quad\varepsilon\not=0.$$  Obviously $ N(A_\varepsilon ) =\{0\},$ and  so $M(A_\varepsilon ) =B(R^2);$
while
$\lim\limits_{\varepsilon\rightarrow 0}A_\varepsilon=A $\ \
and \ \ $\dim M(A_\varepsilon)= 4 $  for
$\varepsilon\not=0.$ Thus  $A_\varepsilon$ for
$\varepsilon\not=0$ is not in $J(A,E_*)$ because  $\dim E_*=1$ and $M(A_\varepsilon ) =B(R^2)$, So that $A$ is
not an inner point of $J(A,\mathbf{E}_*)$.
Therefore $J(A,\mathbf{E}_*)$ is nontrivial.
In this case by Theorem $2.1$  we have to try the co-final set $J(A,\mathbf{E}_*)$ to get  integral submanifold of $\mathcal{F}$
 at $A$.

 {\bf Lemma $4.1$}\quad{\it Suppose that $X\in\Lambda$ is double
  splitting, say that  $X^+$ is a generalized inverse of $X$. Then
  $$M(X)=\left\{P^{N(X^+)}_{R(X)}T+R^{R(X)}_{N(X^+)}TP^{N(X)}_{R(X^+)}:\forall
  T\in B(E ,F),\right \}\eqno(4.1)$$ and one of its complements is $$\mathbf {E}_X=\{P_{N(X^+)}^{R(X)} T P_{N(X)}^{R(X^+)}\}=\mathbf{F}_X .\eqno(4.2) $$  Where $\mathbf{F}_X=\{T\in B(E,F): R(T)\subset N(X^+)$ and $N(T)\supset R(X^+)\}.$

  \textbf{proof}\quad  Obviously $$ T=P^{N(X^+)}_{R(X)}T+P^{R(X)}_{N(X^+)}TP^{N(X)}_{R(X^+)}+ P_{N(X^+)}^{R(X)} T P_{N(X)}^{R(X^+)}\eqno(4.3)$$ for any $T\in B(E,F).$   By $(4.3)$ it is easy to see that $(4.1)$ holds. Especially,   $T\in M(X)$ if and only if $T=P^{N(X^+)}_{R(X)}T+R^{R(X)}_{N(X^+)}TP^{N(X)}_{R(X^+)}$. Indeed, let $T=P_{R(X)}^{N(X^+)}W+P_{N(X^+)}^{R(X)}WP_{R(X^+)}^{N(X)}$ for $W\in B(E,F)$,then $$P_{R(X)}^{N(X^+)}T+P_{N(X^+)}^{R(X)}TP_{R(X^+)}^{N(X)}=P_{R(X)}^{N(X^+)}W+P_{N(X^+)}^{R(X)}WP_{R(X^+)}^{N(X)}=T.$$  Next go to prove $(4.2).$ Let $T$ be any one in $F_X$, then $T$ satisfies  $R(T)\subset N(X^+)$ and $ N(T)\supset R(X^+) ,$ and  by $(4.3)$ $T=P_{N(X^+)}^{R(X)}TP_{N(X)}^{R(X^+)},$ i.e.,
    $E_X \supset F_X.$ The converse relation is immediate from the definition of $E_X$.  Now to end  the proof we  need only to show
  $  M(X)\cap E_X=\{0\}.$ Clearly $M(X)+E_X=B(E,F);$ while $$T= P^{N(X^+)}_{R(X)}T+R^{R(X)}_{N(X^+)}TP^{N(X)}_{R(X^+)}=0 \, \, \mbox{for} \, \, T\in M(X)\cap E_X .$$
 The proof ends. $\Box $

\textbf{Theorem $4.1$}\quad{\it  Suppose that $A\in\Lambda$ is double
splitting, then $\mathcal{F}$ at $A$ is smooth integrable.}

\textbf{Proof}\quad  We will apply Frobenius theorem in Banach spaces
to the proof of the theorem. Hence we have to investigate the
co-final set of $\mathcal{F}$ at $A$ \mbox{and} the coordinate operator
$\alpha(X)$ for $X\in\Lambda$. Let $A^+$ be a generalized inverse of
$A$, and
$\mathbf{E}_*=\{T\in B(E,F):{R}(T)\subset N(A^+)$
and $ N(T)\supset {R}(A^+)\}.$
 By Lemma $4.1$  $M(A)\oplus \mathbf{E}_*=B(E,F).$
 In order to determine $\alpha(X)$ for $X\in J(A,\mathbf{E}_*)$, we claim the following equalities for $X\in J(A,\mathbf{E}_*):$
$$\mathbf{P}^{\mathbf{E}_*}_{M(X)}T=P^{N(A^+)}_{R(X)}T+P^{R(X)}_{N(A^+)}TP^{N(X)}_{R(A^+)}$$
and
$$\mathbf{P}^{M_0}_{\mathbf{E}_*}T=P^{R(A)}_{N(A^+)}TP^{R(A^+)}_{N(A)} \, \, \mbox{for} \, \, T\in B(E,F),$$ for simplicity, hereafter write $M_0=M(A)$.
Let $$\mathbf{P}T=P^{N(A^+)}_{R(X)}T+P^{R(X)}_{N(A^+)}TP^{N(X)}_{R(A^+)} \, \,\forall
T\in B(E,F).$$ Evidently,
\begin{eqnarray*}
\mathbf{P}^2T&=&\mathbf{P}(P^{N(A^+)}_{R(X)}T+P^{R(X)}_{N(A^+)}TP^{N(X)}_{R(A^+)})\\
&=&P^{N(A^+)}_{R(X)}(P^{N(A^+)}_{R(X)}T+P^{R(X)}_{N(A^+)}TP^{N(X)}_{R(A^+)})+P^{R(X)}_{N(A^+)}(P^{N(A^+)}_{R(X)}T+P^{R(X)}_{N(A^+)}TP^{N(X)}_{R(A^+)})P^{N(X)}_{R(A^+)}\\
&=&P^{N(A^+)}_{R(X)}T+P^{R(X)}_{N(A^+)}TP^{N(X)}_{R(A^+)}\\
&=&\mathbf{P}T\quad\quad\forall T\in B(E,F).
\end{eqnarray*}
By Lemma $4.1,$
\begin{eqnarray*}
M(X)&=&\{P^{N(A^+)}_{R(X)}T+P^{R(X)}_{N(A^+)}TP^{N(X)}_{R(A^+)}:\forall
T\in B(E,F)\}\\
&=&\{\mathbf{P}T:\forall T\in B(E,F)\}=R(\mathbf{P});
\end{eqnarray*}
while
\begin{eqnarray*} T\in N(\mathbf{P})&\Leftrightarrow&
P^{N(A^+)}_{R(X)}T+P^{R(X)}_{N(A^+)}TP^{N(X)}_{R(A^+)}=0\\
&\Leftrightarrow&P^{N(A^+)}_{R(X)}T=0 \, \mbox{ and} \, \,
P^{R(N)}_{N(A^+)}TP^{N(X)}_{R(A^+)}=0\\
&\Leftrightarrow&R(T)\subset N(A^+) \, \mbox{ and} \, \, N(T)\supset
R(A^+)\\
&\Leftrightarrow& T\in\mathbf{E}_*. \end{eqnarray*}
This shows $\mathbf{P}=\mathbf{P}^{\mathbf{E}_*}_{M(X)}$. Similarly,
 $\mathbf{P}^{M_0}_{\mathbf{E}_*}T=P^{R(A)}_{N(A^+)}TP^{R(A^+)}_{N(A)}.$
By the same way as that of the proof of the equality
$\alpha(x)=P^{N_0}_{E_*}P^{E_*}_{N_x}\forall x\in U_0$ in Theorem
$3.2,$ one can prove
$$\alpha(X)=\mathbf{P}^{M_0}_{E_*}\mathbf{P}^{E_*}_{M(X)} \, \forall X\in
J(A,\mathbf{E}_*)\eqno(4.5)$$
 too.
In order to get a neighborhood $V$ at $A$ in $M_0$ and a smooth map $\Psi $ from $V$ into $E_*$ such that the conditions (ii) and (iii) in Theorem $2.3 $ hold, we consider such a set $S$ in $J(A,E_*)$ that $\alpha $ for $X\in S$  can be calculated.
Let $S$ be the set of all double splitting operators $X$ in
$W=\{T\in B(E,F):\| T-A\|<\|A^+\|^{-1}\}$ with the generalized
inverse $X^+=A^+C^{-1}_A(A^+,X).$
Clearly $R(X^+)=R(A^+)$ and $N(X^+)=N(A^+)$ because of
$C_A(A^+,X)N(A^+)=N(A^+).$ This is the show of the condition (ii).
Let $E_*=E_A.$ Since $N(X^+)=N(A^+)$ and $R(X^+)=R(A^+)$ for $X \in S,$  it follows $\mathbf{E_X}=\mathbf{E}_*$ and $M(X)\oplus\mathbf{E}_*=B(E,F)$ for $X \in S$ from Lemma $4.1$. Therefore $J(A,\mathbf{E}_*)\supset S.$
For simplicity write
$$P^{R(A)}_{N(A^+)}, P^{N(A^+)}_{R(A)},
 P^{R(A^+)}_{N(A)} \,  \mbox{and}  \, \,  P^{N(A)}_{R(A^+)}$$
 as
 $$P_{N(A^+)},
 P_{R(A)}, P_{N(A)} \,  \mbox{and}  \,  P_{R(A^+)},$$
 respectively, in the sequel. The following equalities for $T\in W $ will be used many times: $$C^{-1}_A(A^+,T)TP_{R(A^+)}=A \, \mbox{and} \,
 C^{-1}_A(A^+,T)P_{N(A^+)}=P_{N(A^+)}\eqno(4.6)$$
 This is immediate from $C_A(A^+,T)A=TP_{R(A^+)} \, \mbox{and} \,
 C_A(A^+,T)P_{N(A^+)}=P_{N(A^+)}.$
Due to the fine property of $S$, we have
$$P^{N(X^+)}_{R(X)}=P^{N(A^+)}_{R(X)}=XA^+C^{-1}_A(A^+,X) \, \mbox{and} \,
 P^{N(X)}_{R(X^+)}=P^{N(X)}_{R(A^+)}=A^+C^{-1}_A(A^+,X)X$$ for $X\in S.$
 Then by $(4.5)$ we  have for each $X\in S$ and all $\Delta X\in M_0$
 \begin{eqnarray*}
 \mathbf{\alpha}(X)\Delta
 X&=&\mathbf{P}^{M_0}_{\mathbf{E}_*}\mathbf{P}^{\mathbf{E}_*}_{M(X)}\Delta
 X\\
 &=&\mathbf{P}^{M_0}_{\mathbf{E}_*}(P^{N(A^+)}_{R(X)}\Delta
 X+P^{R(X)}_{N(A^+)}\Delta XP^{N(X)}_{R(A^+)})\\
 &=&P_{N(A^+)}P^{N(A^+)}_{R(X)}\Delta
 XP_{N_0}+P_{N(A^+)}P^{R(X)}_{N(A^+)}\Delta
 XP^{N(X)}_{R(A^+)}P_{N_0}\\
 &=&P_{N(A^+)}XA^+C^{-1}_A(A^+,X)\Delta XP_{N_0}\\
 &&+P_{N(A^+)}(I_F-XA^+C^{-1}_A(A^+,X))\Delta
 XA^+C^{-1}_A(A^+,X)XP_{N_0}\end{eqnarray*}
 where $N_0=N(A)$.  Meanwhile
 \begin{eqnarray*}
 &&P_{N(A^+)}XA^+C^{-1}_A(A^+,X)\Delta XP_{N_0}\\
 &&\quad =P_{N(A^+)}(P_{N(A^+)}+XA^+)C^{-1}_A(A^+,X)\Delta
 XP_{N_0}\\
 &&\quad\quad-P_{N(A^+)}C^{-1}_A(A^+,X)\Delta XP_{N_0}\\
 &&\quad=P_{N(A^+)}\Delta XP_{N_0}-P_{N(A^+)}C^{-1}_A(A^+,X)\Delta
 XP_{N_0}\\
 &&\quad=-P_{N(A^+)}C^{-1}_A(A^+,X)\Delta XP_{N_0}\end{eqnarray*}
 because $P_{N(A^+)}\Delta XP_{N_0}=0$ for $\Delta X\in M_0$, and
\begin{eqnarray*}
&&P_{N(A^+)}(I_F-XA^+C^{-1}_A(A^+,X))\Delta
XA^+C^{-1}_A(A^+,X)XP_{N_0}\\
&&\quad=P_{N(A^+)}(I_F-(P_{N(A^+)}+XA^+)C^{-1}_A(A^+,X))\Delta
XA^+C^{-1}_A(A^+,X)XP_{N_0}\\
&&\quad\quad+P_{N(A^+)}C^{-1}_A(A^+,X)\Delta
XA^+C^{-1}_A(A^+,X)XP_{N_0}\\
&&\quad=P_{N(A^+)}C^{-1}_A(A^+,X)\Delta
XA^+C^{-1}_A(A^+,X)XP_{N_0}.\end{eqnarray*} Thus
$$\mathbf{\alpha}(X)\Delta X=P_{N(A^+)}(C^{-1}_A(A^+,X)\Delta
XA^+C^{-1}_A(A^+,X)X-C^{-1}_A(A^+,X)\Delta X)P_{N_0}\eqno(4.7)$$ for
$X\in S,$ which is just we please. In order to set the smooth map $\Psi ,$ we consider the smooth map $\mathbf{D}$ from $W$
onto $D(W)$ as follows,
$$T=\mathbf{D}(X)=(X-A)P_{R(A^+)}+C^{-1}_A(A^+,X)X .$$
Obviously, $\mathbf{D}(A)=A$ and
$$\mathbf{D}'(X)\Delta X=\Delta XP_{R(A^+)}+C^{-1}_A(A^+,X)\Delta
X-C^{-1}_A(A^+,X)\Delta XA^+C^{-1}_A(A^+,X)X.$$
Hereby
$$\mathbf{\alpha}(X)\Delta X=P_{N(A^+)}(-\mathbf{D}'(X)\Delta
X)P_{N_0}.\eqno(4.8)$$
In view of $(4.8)$ and $(2.4)$, one can try $\mathbf{D}$ to get the smooth
solution of $(2.4).$
For this we  go to  discuss the map $\mathbf{D}$ more carefully.
Let
$\mathbf{D}_*(X)= XP_{R(A^+)}+C_A(A^+,X)XP_{N_0} \, \forall X\in
B(E,F).$ We claim
$\mathbf{D}_*(X)=\mathbf{D}^{-1}(X) \,\mbox{
for} \, X\in V_1=\{T\in B(E,F):\|(X-A)A^+\|<1\}.$ Clearly
$(\mathbf{D}_*(X)-A)A^+=(X-A)A^+ \, \forall X\in B(E,F).$
Note $AA^+=C^{-1}_A(A^+,X)XA^+$ for $X\in V_1,$ which follows  from
the equality $A=C^{-1}_A(A^+,X)XP_{R(A^+)}$ in $(4.6),$ then one
observes
 \begin{eqnarray*}
 (D(X)-A)A^+&=&(X-A)A^++C^{-1}_A(A^+,X)XA^+-AA^+\\
 &=&(X-A)A^++AA^+-AA^+\\
 &=&(X-A)A^+,\quad\quad\forall X\in V_1.
 \end{eqnarray*}
Hereby
$$C_A(A^+,\mathbf{D}(T))=I_F+(D(T)-A)A^+=C_A(A^+,T)$$
and
$$C_A(A^+,\mathbf{D}_*(T))=I_F+(\mathbf{D}_*(T)-A)A^+=C_A(A^+,T)$$
for all $T\in V_1.$
Consequently  the following equalities for all $T \in V$ hold :
$$  C^{-1}_A(A^+,\mathbf{D}(T))=C^{-1}_A(A^+,T) \, \, \mbox{and} \, \,
C^{-1}_A(A^+,\mathbf{D}_*(T))= C^{-1}_A(A^+,T).\quad \eqno (4.9) $$
 Then by the equality $C^{-1}_A(A^+,X)XP_{R(A^+)}=A$ in $(4.6),$
 \begin{eqnarray*}
 (\mathbf{D}_*\circ\mathbf{D})(X)&=&\mathbf{D}(X)P_{R(A^+)}+C_A(A^+,\mathbf{X}))\mathbf{D}(X)P_{N_0}\\
 &=&(X-A)P_{R(A^+)}+C^{-1}_A(A^+,X)XP_{R(A^+)}\\
 &&+C_A(A^+,X)\mathbf{D}(X)P_{N_0}\\
 &=&XP_{R(A^+)}-A+A+XP_{N_0}=X,\quad\forall X\in V_1.\end{eqnarray*}
 Similarly
 $$(\mathbf{D}\circ\mathbf{D}_*)(T)=T,\quad\forall T\in V_1.$$
 Thus  $\mathbf{D}_*(X)=\mathbf{D}^{-1}(X)\ \forall X\in V_1$ and $\mathbf{D}$ is a smooth diffeomorphism from $V_1$ onto itself
 with $\mathbf{D}(A)=A.$ (The diffeomorphisn $D$ is first presented in [Ma$5$].)
Let
$V_*=\mathbf{D}(W)\cap M_0  \, \mbox{and}  \,
\mathbf{\Phi}(T)=\mathbf{D}_*(T) \, \forall T\in V_*.$ Clearly
$V_*$ is a neighborhood at $A$ in $M_0$, and  $\mathbf{\Phi}(T)$ for
$T\in V_*$ is a $C^\infty$ map.
Directly
 \begin{eqnarray*}
 C^{-1}_A(A^+,\mathbf{\Phi}(T))\mathbf{\Phi}(T)&=&C^{-1}_A(A^+,T)\mathbf{\Phi}(T)\\
 &=&C^{-1}_A(A^+,T)TP_{R(A^+)}+TP_{N_0}\\
 &=&A+TP_{N_0},\quad\quad\forall T\in V_*.\end{eqnarray*}
Note that $R(TP_{N_0})\subset R(A)$ because of $T\in V_*\subset
M_0.$ By the equivalence of the conditions (ii) and (vii) in Theorem
$ 1.1,$ one can conclude
 ${X}=\mathbf{\Phi}(T)\in S \, \mbox{for} \, T\in V_*.$
Moreover by $(4.7)$ and $(4.9)$
 $$\begin{array}{rllr}
 \alpha(\mathbf{\Phi}(T))\Delta X&=&P_{N(A^+)}(C^{-1}_A(A^+,T)\Delta
 XA^+C^{-1}_A(A^+,T)\mathbf{\Phi}(T)\\
 &&-C^{-1}_A(A^+,T)\Delta X)P_{N_0}\\
 &=&P_{N(A^+)}C^{-1}_A(A^+,T)\Delta XA^+TP_{N_0}\\
 &&-P_{N(A^+)}C^{-1}_A(A^+,T)\Delta XP_{N_0},\quad\quad\forall T\in
 V_*.\end{array}\eqno(4.10)$$
 Let
 $$\mathbf{
 \Phi}_0(T)=\mathbf{P}^{\mathbf{E}_*}_{M_0}\mathbf{\Phi}(T)
 \quad \mbox{and}\quad
 \mathbf{\Phi}_1(T)=\mathbf{P}^{M_0}_{\mathbf{E}_*}\mathbf{\Phi}(T).$$
 Since $\Phi$ is smooth, both $\Phi _0, \Phi _1 $ are smooth too .Then we have $\mathbf{\Phi}_0(A)=A$ and
 \begin{eqnarray*}
 \mathbf{\Phi}_0(T)&=&P_{R(A)}\mathbf{\Phi}(T)+P_{N(A^+)}\mathbf{\Phi}(T)P_{R(A^+)}\\
 &=&TP_{R(A^+)}+P_{R(A)}C_A(A^+,T)TP_{N_0}\end{eqnarray*}
 for $T\in V_1.$ Similarly,
 $$\mathbf{\Phi}_1(T)=P_{N(A^+)}C_A(A^+,T)TP_{N_0}\quad\quad\forall
 T\in V_1.$$
 By direct computing,
 $$\mathbf{\Phi}'_0(T)\Delta T=\Delta TP_{R(A^+)}+P_{R(A)}\Delta
 TA^+TP_{N_0}+P_{R(A)}TA^+\Delta TP_{N_0}$$
 and
$$\begin{array}{rll}\mathbf{\Phi}'_1(T)\Delta T&=&P_{N(A^+)}\Delta
TA^+TP_{N_0}+P_{N(A^+)}C_A(A^+,T)\Delta TP_{N_0},\\
&=&P_{N(A^+)}(\Delta TA^+T+TA^+\Delta T)P_{N_0}\quad\quad\forall
T\in V_1.\end{array}$$ (Note $P_{N(A^+)}\Delta TP_{N_0}=0$ for
$\Delta T\in M_0$.)
 Consider the map $\mathbf{\Phi}_0|_{V_*},$ still
write as $\mathbf{\Phi}_0$ for simplicity. Obviously,
$$\mathbf{\Phi}'_0(A)\Delta T=\Delta TP_{R(A^+)}+\Delta TP_{N_0}=\Delta
T,\quad\quad \forall \Delta T\in M_0.$$ Then by
 the inverse map theorem, there exists a neighborhood $V_0$ at $A$
 such that
$\mathbf{\Phi}_0$ from $V_0$ onto $\mathbf{\Phi}_0(V_0)$ is a smooth
diffeomorphism. We are now in the position to set the neighborhood $V$ and the smooth map $\Psi $ in $(2.4).$ Let
$$V=\mathbf{\Phi}_0(V_0) \,  \mbox{and} \,
\mathbf{\Psi}(Z)=(\mathbf{\Phi}_1\circ\mathbf{\Phi}^{-1}_0)(Z) \, \mbox{
for} \,  Z\in V.$$
It follows that $\Psi $ is smooth from the smooth diffeomorphism $\Phi_0$ from $V_0$ onto $\Phi_0(V_0).$  Finally we are going to prove that
 $\mathbf{\Psi}(T)$ for $T\in V$ is the solutions of $(2.4). $ Take the transformation
 $Z=\mathbf{\Phi}_0(T)$ for any $T\in V$, then  we have
 $$\Delta T=(\mathbf{\Phi}'_0(T))^{-1}\Delta Z=(\mathbf{\Phi}^{-1}_0)'(Z)\Delta Z,\quad\quad\forall\Delta Z\in M_0.$$
It is indicated in the above that $\mathbf{\Phi}(T)\in S \,\mbox{ for \ any}  \, T\in V_*.$
Hence
$Z+\mathbf{\Psi}(Z)=\mathbf{\Phi}_0(T)+\mathbf{\Phi}_1(T)=\mathbf{\Phi}(T)\in
S \, \forall Z\in V \, \mbox{because \,of } V_0\subset V_*$. So
   $$\mathbf{\alpha}(Z+\mathbf{\Psi}(Z))=\mathbf{\alpha}(\mathbf{\Phi}_0(T)+\mathbf{\Phi}_1(T))=\mathbf{\alpha}(\mathbf{\Phi}(T)), \forall Z\in V.$$
Note $C^{-1}_A(A^+,T)P_{N(A^+)}=P_{N(A^+)}$ and
$P_{R(A)}=I_F-P_{N(A^+)}.$ By $(4.10)$
 $$\begin{array}{rllr}
&&\alpha(\mathbf{\Phi}(T))\mathbf{\Phi}'_0(T)\Delta T\\
&&\quad=P_{N(A^+)}C^{-1}_A(A^+,T)(\mathbf{\Phi}'_0(T)\Delta T)A^+TP_{N_0}\\
&&\quad\quad-P_{N(A^+)}C^{-1}_A(A^+,T)(\mathbf{\Phi}'_0(T)\Delta
T)P_{N_0}\\
&&\quad=P_{N(A^+)}C^{-1}_A(A^+,T)\Delta TA^+TP_{N_0}\\
&&\quad\quad-P_{N(A^+)}C^{-1}_A(A^+,T)P_{R(A)}(\Delta
TA^+T+TA^+\Delta T)P_{N_0}\\
&&\quad=P_{N(A^+)}C^{-1}_A(A^+,T)P_{N(A^+)}(\Delta
T)A^+TP_{N_0}\\
&&\quad\quad+P_{N(A^+)}TA^+(\Delta T)P_{N_0}.
\end{array}$$
Then by  the equality in $(4.6),$
$$C^{-1}_A(A^+,T)P_{N(A^+)}=P_{N(A^+)},$$
we  further have $$\begin{array}{rl}
&\alpha(\mathbf{\Phi}(T))\mathbf{\Phi}'_0(T)\Delta T\\
&\quad=P_{N(A^+)}\{(\Delta T)A^+T+TA^+(\Delta T)\}P_{N_0}\\
&\quad=\mathbf{\Phi}'_1(T)\Delta T. \end{array} $$ Therefore
$$\alpha(Z+\mathbf{\Psi}(Z))\Delta
Z=\alpha(\mathbf{\Phi}(T))\mathbf{\Phi}'_0(T)\Delta
T=\mathbf{\Phi}'_1(T)\Delta T.$$ Now take
$$\Delta T=(\mathbf{\Phi}^{-1}_0)'(Z)\Delta Z.$$
Then one observes
$$\alpha(Z+\mathbf{\Psi}(Z))=\mathbf{\Psi}'(Z),\quad\forall Z\in
V,$$ which  also implies that the condition (ii) in Theorem $2.3$ holds since $\mathbf{\Psi}$ is a smooth map.
Obviously
\begin{eqnarray*}
\mathbf{\Psi}(\mathbf{P}^{E_*}_{M_0}A)&=&\mathbf{\Psi}(A)\\
&=&\mathbf{\Phi}_1(\mathbf{\Phi}_0(A))\\
&=&\mathbf{\Phi}_1(A)=\mathbf{P}^{M_0}_{\mathbf{E}_*}A=0.
\end{eqnarray*}
So far the theorem is proved. $\Box$

 Specially we have the global result as follows,

{\bf Theorem 4.2}\quad{\it  Each of $F_k,\Phi_{m,n},\Phi_{m,\infty}$
  and $\Phi_{\infty,n}$ is a smooth submanifold in $B(E,F)$ and tangent to $M(X)$ at any $X$ in it.}

{\bf Proof}\quad Let $\mathbf{F}$ be any one of
  $F_k,\Phi_{m,n},\Phi_{m,\infty}$ and $\Phi_{\infty,n},$ and $A$
  any one operator in $\mathbf{F}$.
It is clear that $A$ is double splitting, say that $A^+$ is a
generalized inverse of $A$. Write
$$\begin{array}{rl}
&V=\{T\in B(E,F):\|(T-A)A^+\|<1\},\\
&W=\{T\in B(E,F):\|T-A\|<\|A^+\|^{-1}\}\end{array}$$
and
$M_0=\{T\in B(E,F):TN(A)\subset R(A)\}.$
It is proved in the proof of Theorem $4.1$ that
$M_0\oplus\mathbf{E}_*=B(E,F)$ where
$\mathbf{E}_*=\{T\in B(E,F):R(T)\subset N(A^+) \, \mbox{and} \, N(T)\supset R(A^+)\}.$
Clearly
$\mathbf{D}(X)=(X-A)P_{R(A^+)}+C^{-1}_A(A^+,X)X $
 for $X\in V$ is a $C^\infty$ diffeomorphism from $V$ onto itself with
$\mathbf{D}(A)=A$ and
$\mathbf{D}^{-1}(T)=TP_{R(A^+)}+C_A(A^+,T)TP_{N(A)}.$
We claim
$$\mathbf{D}(\mathbf{F}\cap W)=M_0\cap \mathbf{D}(W).$$
By Theorem $1.2$ and the equivalence of the conditions (i) and (vi) in
Theorem$1.1,$
$$\mathbf{D}(X)N(A)=C^{-1}_A(A^+,X)X{N(A)}\subset R(A)$$
for any $X\in\mathbf{F}\cap W.$
So
$\mathbf{D}(\mathbf{F}\cap W)\subset M_0\cap\mathbf{D}(W).$
Conversely
$$\begin{array}{rllr}
&&C^{-1}_A(A^+,\mathbf{D}^{-1}(T))\mathbf{D}^{-1}(T)N(A)\\
&&\quad=C^{-1}_A(A^+,T)\mathbf{D}^{-1}(T)N(A)\\
&&\quad=C^{-1}_A(A^+,T)C_A(A^+,T)TN(A)\\
&&\quad=TN(A)\subset R(A)\quad\forall T\in M_0\cap\mathbf{D}(W).
\end{array}$$
Then by Theorem $1.2$ and the equivalence of the conditions (i) and (vi) in
Theorem $1.1,$ we conclude
$X=\mathbf{D}^{-1}(T)\in\mathbf{F}\cap W \,\forall T\in
M_0\cap\mathbf{D}(W).$
So
$\mathbf{D}(X)=T \,\forall T\in M_0\cap\mathbf{D}(W),$
i.e.,
$\mathbf{D}(\mathbf{F}\cap W)\supset M_0\cap\mathbf{D}(W).$
This proves
$\mathbf{D}(\mathbf{F}\cap W)=M_0\cap \mathbf{D}(W).$

We now see that $(\mathbf{D}, W,B(E,F))$ at $A$ is a local
coordinate chart of $B(E,F)$ satisfying that
$\mathbf{D}(\mathbf{F}\cap W)=M_0\cap\mathbf{D}(W)$ is an open set
in $M_0$, and $M_0\oplus\mathbf{E}_*=B(E,F)$. So $\mathbf{F}$ is a
submanifold in $B(E,F)$.
Next to show that $\mathbf{F}$ is smooth submanifold in \(B(E,F)\). Let
$W_1=\{T\in B(E,F):\|T-B\|<\|B^+\|^{-1}\},$
and
 $M_1=\{T\in B(E,F):TN(B)\subset R(B)\}.$
 Obviously
 $\mathbf{D}_1(T)=(T-B)P_{R(B^+)}+C^{-1}_B(B^+,T)T$
 is a $c^\infty$ diffeomorphism from $W_1$ onto $\mathbf{D}_1(W_1)$
 with $\mathbf{D}_1(B)=B$. When $\mathbf{F}\cap W\cap
 W_1\not=\varnothing,$  it is clear that $\mathbf{D}\circ
 \mathbf{D}^{-1}_1$ from $\mathbf{D}_1(\mathbf{F}\cap W\cap W_1)$ onto
 $\mathbf{D}(\mathbf{F}\cap W\cap W_1)$ is of $C^\infty$. So  we
 conclude $\mathbf{F}$ is smooth submanifold in $B(E,F)$.
In the same way as the proof of the conclusion \(\mathbf{D}(\mathbf{F}\cap W)=M_0\cap \mathbf{D}(W)\) in the above, we also infer
$\mathbf{D}_1(\mathbf{F}\cap W_1)=M_1\cap\mathbf{D}_1(W_1),$
and $\mathbf{D}_1$ from $W_1$ onto $\mathbf{D}_1(W_1)$ is smooth
$C^\infty$ diffeomorphism with $\mathbf{D}_1(B)=B$. Evidently
$\mathbf{D}_1(\mathbf{F}\cap W\cap W_1)=M_1\cap\mathbf{D}_1(W\cap
W_1),$ and
$\mathbf{D}(\mathbf{F}\cap W\cap W_1)=M_0\cap\mathbf{D}(W\cap
W_1),$ so that
$\mathbf{D}\circ\mathbf{D}^{-1}_1 \, \mbox{ from} \,
M_1\cap\mathbf{D}_1(W\cap W_1) \, \mbox{ onto} \,
M_0\cap\mathbf{D}(W\cap W_1)$ is smooth diffeomorphism when
$\mathbf{F}\cap W\cap W_1\not=\varnothing$. This proves that
$\mathbf{F}$ is smooth submanifold in $B(E,F)$. Finally go to show
that $\mathbf{F}$ is tangent to $M(X)$ at any $X$ in $\mathbf{F}$.
Since $A$ is an arbitrary operator in $\mathbf{F}$,it is enough to show that $\mathbf{F}$ at $A$ is tangent to $M(A)$. Let $T(X)=\{ \dot{c(0)}:
c^1-curve \, \,c(t)\subset S \, \mbox{with} \, \, c(0)=X.\} $ It is proved in the above that for any
$A\in\mathbf{F}$ there are the neighborhood $W$ at $A$, the space
 $\mathbf{E}_*$, and the diffeomorphism $\mathbf{D}$ from
$W$ onto $\mathbf{D}(W)$ satisfying
$M(A)\oplus\mathbf{E}_*=B(E),\,\mathbf{D}(\mathbf{F}\cap
W)=M(A)\cap\mathbf{D}(W) \, \mbox{and} \,
\mathbf{D}'(A)=\mathbf{I}.$ By the similar way to the proof of the
equality $\varphi'(x)M(x)=E_0\forall x\in S\cap U_0$ in Theorem $2.1,$
one can prove
$\mathbf{D}'(X)T(X)=M(A), \, \forall X\in\mathbf{F}\cap W$
and so $T(A)=M(A).$ This means that $\mathbf {F}$ at $A$ is tangent to $M(A).$  Then proof ends. \quad$\Box$.

These results seem to be meaningful for differential topology, global analysis and the geometrical method in
differential equations (see [Abr],[An],[Caf] and [Bo]).\\ \quad

Let $U$ be an open set in E, $S$ a $c^1$ submanifold in $U$, and $f$ a nonlinear functional from $U$ into $(-\infty ,\infty )$. By Theorem $2.2$ there exist a neighborhood $U_0$ at $x_0$ in $E$, a subspace $E_*$ and a diffeomorphism $\varphi $ from $U_0$ onto $\varphi (U_0),$ such that $\varphi (S\cap U_0)=V_0$ is an open set in $T_{x_0}S,$ $ T_{x_0}S\oplus E_*=E,$ $ \varphi'(x_0)=I_E$  and $ \varphi '(x) T_xS= T_{x_0}S \, \mbox{for} \, x\in S\cap U_0.$ Naturally, a point $x_0$ is said to be critical point of $f$ under the constraint $S$ provided $\varphi(x_0)$ is the point of $(f \circ \varphi ^{-1})(x): V_0\rightarrow (-\infty,\infty).$ certainly $(f\circ \varphi ^{-1})'(\varphi (x_0))=0$ in $ B(T_{x_0}S, (-\infty ,\infty )),$ i.e., $$(f\circ\varphi^{-1})'(\varphi(x_0))e=0 \quad \forall e\in T_{x_0}S.$$ Note $\varphi'(x_0)=I_E.$ Directly  $$f'(x_0)(\varphi^{-1})'(\varphi(x_0))e=f'(x_0)\varphi'(x_0)^{-1}e=f'(x_0)e=0 \, \forall e\in T_{x_0}S.$$ So $N(f'(x_0))\supset T_{x_0}S. $ Therefore we have\\

 \textbf {Theorem $4.3$}\quad{\it If $x_0\in U$ is a critical point of $f$ under the constraint of $c^1$ Banach submanifold $S$, then $$ N(f'(x_0))\supset T_{x_0}S.$$ } \ The theorem expands the principle for critical point of $f$ under the generalized regular constraint to the constraint of $c^1$ submanifold in Banach spaces $S$. The examples $1-4$ in [Ma $8$] illustrate the approach with Theorem $4.3$ to get critical point of $f$ under the constraint of $c^1$ Banach submanifold $S$. In view of the above Banach submanifolds with the explicit expression of the tangent space ,the principle should be potential in applcation.

\newpage
\begin{center}\textbf{Appendix}
\end{center}

\textbf{$1$ \quad The proof of Theorem $1.1$}\\ Go to show $(vi)\Rightarrow(ii).$ M. Z. Nashed and
Chen have proved $(vi) \Rightarrow (ii)$ in [N-C]. In fact, the inverse relation hold too.
 Note $C^{-1}(A^+,T)TA^+A=A \forall T \in T V(A,A^+)$.
Then we have $$BTB-B=0 \, \,\mbox{and} \, \,
TBT-T=-(I_F-AA^+)C^{-1}(A^+,T)T\eqno (1)$$ for any $T\in V(A,A^+)$.
So(vi)$\Leftrightarrow$(ii).

Evidently,
$$\begin{array}{rllr}
C^{-1}_A(A^+,T)Th&=&C^{-1}_A(A^+,T)TA^+Ah+C^{-1}_A(A^+,T)T(I_E-A^+A)h\\
&=&Ah+C^{-1}_A(A^+,T)T(I_E-A^+A)h,\quad\forall h\in E,\end{array}$$
so that (vi)$\Leftrightarrow$(vii).

Go to the show (v)$\Leftrightarrow$(vi). Assume that (vi) holds.
We have for any $h\in N(A)$, there exists $g\in R(A^+)$ such that
$C^{-1}_A(A^+,T)Th=Ag$. Note $$A^+Ag=g \mbox{ and } C^{-1}_A(A^+,T)TA^+A=A,$$ then
$$C^{-1}_A(A^+,T)Th=Ag=C^{-1}_A(A^+,T)TA^+Ag=C^{-1}_A(A^+,T)Tg.$$ So
$h-g\in N(T)$ and satisfies $(I_E-A^+A)(h-g)=h$. This shows that (v)
holds because of $(I_E-A^+A)N(T)\subset N(A).$ Conversely, assume that (v) holds. Then we have for any $h\in
N(A)$, there exists $g\in N(T)$ such that $h=(I_E-A^+A)g,$ and so,
$C^{-1}_A(A^+,T)Th=-C^{-1}_A(A^+,T)TA^+Ag=-Ag\in R(A)$, i.e., (vi)
holds. This shows (v)$\Leftrightarrow$ (vi).

Go to the show(i)$\Leftrightarrow$(ii). Obviously,
(ii)$\Rightarrow$(i). In fact $\{0\}=R(T)\cap N(B)=R(T)\cap N(A^+).$
Conversely, assume that (i) holds. Obviously $R(TBT-T)\subset R(T)$,
while by (1) that contains in $N(A^+)$.  Then by (i),
$R(T)\cap N(A^+)=\{0\}$, one concludes $TBT-T=0$. Thus $B$ is the
generalized inverse of $T$, this shows (i)$\Rightarrow$(ii). So
(i)$\Leftrightarrow$(ii).

Go to show (i)$\Leftrightarrow$(iii). Obviously,
(iii)$\Rightarrow$(i). Conversely, assume that (i) holds. Then $B$
is the generalized inverse of $T$ with $N(B)=N(A^+)$ because of the
equivalence of (i) and (ii). So (i)$\Leftrightarrow$(iii).

Finally
 go to show (i)$\Leftrightarrow$(iv). Assume that (i) holds.  Since
 (i)$\Rightarrow$(ii) we have $E=N(T)\oplus R(B)=N(T)\oplus R(A^+),$ i.e.,  (iv) holds.
 Conversely, assume that (iv) holds. Then
 $N(A)=(I_E-A^+A)E=(I_E-A^+A)N(T)$ i.e., (v) holds. Hence (i) holds
 since  (v)
 $\Leftrightarrow$(vi)$\Leftrightarrow$(ii)$\Leftrightarrow$(i).
So (i)$\Leftrightarrow$(iv).
 So far we have proved the following relations
 (i)$\Leftrightarrow$(ii),(i)$\Leftrightarrow$(iii),(i)$\Leftrightarrow$(iv),
 (ii)$\Leftrightarrow$(vi),
 (v)$\Leftrightarrow$(vi), and (vi)$\Leftrightarrow$(vii). Thus the
 proof of the theorem is completed.

\textbf{$ 2$ \quad The proof of Theorem $1.2$}
 Assume that the condition $R(T)\cap
 N(A^+)=\{0\}$ holds for $T\in V(A,A^+)$. By (ii) in Theorem 1.1,
 $T$ has a generalized inverse $B$ with $N(B)=N(A^+)$ and
$R(B)=R(A^+)$,
 so that
 $$N(T)\oplus R(A^+)=E=N(A)\oplus R(A^+)$$
 and
$$R(T)\oplus N(A^+)=F=R(A)\oplus N(A^+).$$
 Hereby one observes that $A$ and $T$
 belong the same class.
 Conversely, assume that $T\in V(A,A^+)$  belong to  any one of
 $F_k,\Phi_{m,n},\Phi_{m,\infty}$ and $\Phi_{\infty,n}$.
By (1)
 $B=A^+C^{-1}_A(A^+,T)=D^{-1}_A(A^+,T)A^+$ satisfies
$BTB=B.$  Thus $B$ and $T$ bear two projections $P_1=BT$ and
$P_2=TB$.
 Indeed, $P^2_1=BTBT=BT=P_1,$ and $P^2_2=TBTB=TB=P_2$.
It  is clear that $T$ is double splitting, say that $T^+$ is a
generalized inverse of $T$. Since $N(B)=N(A^+),$
$$N(P_1)=N(T)\oplus\{ e\in R(T^+): Te\in N(A^+)\}.$$
We next claim
 $$R(P_1)=R(A^+), R(P_2)=R(TA^+),\quad{\rm and}\quad
 N(P_2)=N(A^+).$$

Obviously, $R(P_1)\subset R(B)= R(A^+)$. On the  other hand, by
$C^{-1}_A(A^+,T)TA^+=AA^+$
 $$P_1A^+e=A^+C^{-1}_A(A^+,T)TA^+e=A^+AA^+e=A^+e,\quad \forall e\in
 F.$$
 So
 $R(P_1)=R(A^+)$.

 Obviously,
 $$R(P_2)=R(TA^+C^{-1}_A(A^+,T))=R(TA^+),$$
 $$N(P_2)=N(TB)\supset N(B)= N(A^+),$$
and
 $$N(A^+)=N(B)=N(BTB)=N(BP_2)\supset N(P_2).$$
 So $N(P_2)=N(A^+)$.
 Thus we have
 $$F=R(P_2)\oplus N(P_2)=R(TA^+)\oplus N(A^+)\eqno(2)$$
 and
 $$E=R(P_1)\oplus N(P_1)=R(A^+)\oplus E_*\oplus N(T),\eqno(3)$$
 where
 $E_*=\{e\in R(T^+):Te\in N(A^+)\}.$
 We are now in the position to end the proof.
 Assume that $T\in V(A,A^+)=\{T\in B(E,F):\|T-A\|<\|A^+\|^{-1}\}$ satisfies rank$T=$rank$A<\infty$.
 By $(2)$
 $$F=R(TA^+)\oplus N(A^+)
 =R(A)\oplus N(A^+), \quad \eqno(4)$$
 so that
 $$\dim R(TA^+)(=\dim R(A))=\dim R(T)<\infty.$$
 By $(3)$ and $(4)$
 $$R(T)=R(TA^+)\oplus TE_*$$
 because of $TE_*\subset N(A^+)$. Thus since
 $$\dim R(T)=\dim R(TA^+)<\infty,$$
 $\dim (TE_*)=0$. This shows $R(T)\cap N(A^+)=\{0\}.$
 Assume that $T\in V(A,A^+)$ satisfies $\mathrm{codim} R(T)=\mathrm{codim} R(A)<\infty$, and  both $N(T), N(A)$ are splitting in $E$. Go
 to show $R(T)\cap N(A^+)=\{0\}$. By $(2)$ and $(3)$
 $$R(TA^+)\oplus TE_*\oplus N(T^+)=F=
R(TA^+)\oplus N(A^+).$$
Indeed, by $R(T)=R(TA^+)\oplus TE_*,$ and so $F=R(T)\oplus N(T^+)=R(TA^+)\oplus TE_*\oplus N(T^+).$
Then  $\dim (TE_*\oplus N(T^+))=\dim N(A^+)<\infty$, and so, by the
assumption $\dim N(T^+)=\dim N(A^+)<\infty$, $\dim T(E_*)=0.$ Therefore
$R(T)\cap N(A^+)=\{0\}$.
 Assume that $T\in V(A,A^+)$ satisfies $\dim N(T)=\dim N(A)<\infty$, and both  $R(T), R(A)$ and   splitting in $F$.
  According to $(3),$ we have
 $$R(A^+)\oplus E_*\oplus N(T)=E=R(A^+)\oplus N(A).$$
 Hence $\dim (E_*\oplus N(T))=\dim N(A)<\infty$. Then by the assumption
 $\dim N(T)=\dim N(A)<\infty$, $\dim E_*=0$, so that $R(T)\cap N(A^+)=\{0\}$. Now Theorem $1.2$ is proved.

 \textbf{$3$\quad The Proof of Theorem $1.3$ }

Let $T^+_0$ and $T^\oplus_0$ be two generalized inverses of
$T_0=T_{x_0}$,$\delta=\min\{\|T^+_0\|^{-1},\|T^+_0T_0T^\oplus_0\|^{-1}\|\}$
and $V_\delta=\{T\in B(E,F):\|T-T_0\|<\delta\}.$

Assume for neighborhood $U_0$ at $x_0, R(T_x)\cap
N(T^+_0)=\{0\}\forall x\in U_0$. Since $T_x$ is continuous at $x_0$,
there exists a neighborhood $U_1$ at $x_0$ such that $T_x\in
V_\delta\forall x\in U_1$, then the neighborhood $U_1\cap U_0$ at
$x_0$, for simplicity, still write it as $U_0$, satisfies that
$$T_x\in V_\delta\quad{\rm and}\quad R(T_x)\cap
N(T^+_0)=\{0\}\quad\forall x\in U_0.$$ Write $B=T^+_0T_0T^\oplus_0.$
It is immediate that
$$ BT_0B=B\quad{\rm and}\quad T_0BT_0=T_0,$$
so that $B$ is also a generalized inverse of $T_0$ with
$$N(B)=N(T^\oplus_0)\quad{\rm and}\quad R(B)=R(T^+_0).$$
Indeed
$$f\in N(B)\Rightarrow T_0Bf=T_0T^+_0T_0T^\oplus_0f=T_0T^\oplus
_0f=P^{N(T^\oplus_0)}_{R(T_0)}f=0\Rightarrow f\in N(T^\oplus_0),$$
and
$$e\in R(B)\Leftrightarrow
e=P^{N(T_0)}_{R(B)}e=BT_0e=T^+_0T_0T^\oplus_0T_0e=T^+_0T_0e=P^{N(T_0)}_{R(T^+_0)}e\Leftrightarrow
e\in R(T^+_0).$$ Thus, instead of $A$ and $A^+$ in Theorem $1.1$ by
$T_0$ and $T^+_0$, respectively, the equivalence of the conditions
(i) and (iv) in the theorem shows
$$R(B)\oplus N(T_x)(=R(T^+_0)\oplus N(T_x))=E\quad\forall x\in
U_0.$$ Similarly, instead of $A$ and $A^+$ in Theorem $ 1.1$ by $T_0$
and $B$, respectively, one concludes
$$R(T_x)\cap N(T^\oplus _0)(=R(T_x)\cap N(B))=\{0\}\quad\forall x\in
U_0.$$ The proof ends.

\textbf{$4$ \quad The Proof of Theorem $1.4$}

Assume that $x_0$ is a locally fine point of $T_x$.
Let $T^+_0$ is an arbitrary generalized inverse of $T_0=T_{x_0}.$
 Go to show for $T^+_0$ there is a neighborhood $U_0$ at $x_0$ such that $T_x$ for any $x\in U_0$ has a generalized inverse $T^+_x$ satisfying $\lim\limits_{x\rightarrow x_0}T^+_x=T^+_0.$
  By Definition 1.1, for $T^+_0$ there exists a neighborhood $U_0$ at $x_0$ such that $R(T_x)\cap N(T^+_0)=\{0\}\ \forall x\in U_0.$
  Since $T_x$ is continuous at $x_0$, we can  assume $U_0\subset\{x\in
X:\|T_x-T_0\|<\|T^+_0\|^{-1}\}$. Thus by the equivalence of the
conditions (i) and (ii) in Theorem 1.1,
$T^+_x=T^+_0C^{-1}_{T_0}(T^+_0,T_x)$ is a generalized inverse of
$T_x$ for all $x\in U_0$, and clearly $\lim\limits_{x\rightarrow
x_0}T^+_x=T^+_0$. Conversely, assume that there is a neighborhood
$U_0$ at $x_0$ such that there exists a generalized inverse $T^+_x$
of $T_x$ for all $x\in U_0$, and $\lim\limits_{x\rightarrow
x_0}T^+_x=T^+_0$. Consider  the following projections:
$$P_x=I_E-T^+_xT_x\quad{\rm for}\quad x\in U_0$$
 and $P_0=I_E-T^+_0T_0.$ Obviously $R(P_x)=N(T_x)$ and
 $R(P_0)=N(T_0)$ Let $V_0=\{x\in U_0:\|P_x-P_0\|<1\}\cap\{x\in
 U_0:\|T_x-T_0\|<\|T^+_0\|^{-1}\}.$ Since both $T_x,P_x$ are
 continuous at $x_0,V_0$ is a neighborhood at $x_0$.
 Due to $\|P_x-P_0\|<1 \forall x\in V_0,P_0N(T_x)=N(T_0)$ (see [Ka])
 i.e., $(I_E-T^+_0T_0)N(T_x)=N(T_0)$. Then by the equivalence of the
 conditions (i) and (v),
 $$R(T_x)\cap N(T^+_0)=\{0\},\quad\quad\forall x\in V_0.$$

\textbf{$5$\quad The Proof of Theorem $1.5$}

First go to show that $\alpha$ is unique for which $E_1=\{e+\alpha
e:\forall e\in E_0\}$. If $\{e+\alpha e:\forall e\in
E_0\}=\{e+\alpha_1e:\forall e\in E_0\}$, then for any $e\in E_0$
there exists $e_1\in E_0$ such that $(e-e_1)+(\alpha e-\alpha_1
e_1)=0$, and so, $e=e_1$  and $\alpha e=\alpha_1e, $ i.e.,
$\alpha_1=\alpha$.  This says that $\alpha$ is unique. We  claim
that $\alpha=\left.P^{E_0}_{E_*}P^{E_*}_{E_1}\right|_{E_0}\in
B(E_0,E_*)$
 fulfils $E_1=\{e+\alpha e:\forall e\in E_0\}$.

Obviously,
$$P^{E_*}_{E_0}P^{E_*}_{E_1}e=P^{E_*}_{E_0}(P^{E_*}_{E_1}e+P^{E_1}_{E_*}e)=P^{E_*}_{E_0}e=e,\quad\quad\forall
e\in E_0,$$ and
$$P^{E_*}_{E_1}P^{E_*}_{E_0}e=P^{E_*}_{E_1}(P^{E_*}_{E_0}e+P^{E_0}_{E_*}e)=P^{E_*}_{E_1}e=e,\quad\quad\forall
e\in E_1.$$ Hereby let $\alpha=P^{E_0}_{E_*}P^{E_*}_{E_1}|_{E_0}$,
then
$$e+\alpha e=P^{E_*}_{E_0}P^{E_*}_{E_1}e+P^{E_0}_{E_*}P^{E_*}_{E_1}e=P^{E_*}_{E_1}e\in
E_1$$ for any $e\in E_0$; conversely,
$$e=P^{E_*}_{E_0}e+P^{E_0}_{E_*}e=P^{E_*}_{E_0}e+P^{E_0}_{E_*}P^{E_*}_{E_1}P^{E_*}_{E_0}e$$
for any $e\in E_1.$ So
$$E_1=\{e+\alpha e:\forall e\in E_0\}.$$
Now assume $E_1=\{e_0+\alpha e_0:\forall e_0\in E_0\}$ for any
$\alpha\in B(E_0,E_*)$. Go to show $E_1\oplus E_*=E$. We first claim
that $E_1$ is closed. Let $e_n+\alpha e_n\rightarrow e_*$ as
$n\rightarrow\infty$ where $e_n\in E_0, n=1,2,3,\cdots.$ Then
$$P^{E_*}_{E_0}(e_n+\alpha e_n)=e_n\rightarrow P^{E_*}_{E_0}e_*\in E_0,\mbox{ and}\
 P^{E_0}_{E_*}(e_n+\alpha e_n)=\alpha
e_n\rightarrow\alpha(P^{E_*}_{E_0}e_*)$$ since $\alpha \in B(E_0,E_*)$ and
$P^{E_*}_{E_0}\in B(E)$.  So
$$e_*=P^{E_*}_{E_0}e_*+P^{E_0}_{E_*}e_*=P^{E_*}_{E_0}e_*+\alpha(P^{E_*}_{E_0}e_*)\in
E_1.$$ This shows that $E_1$ is closed. It follows $E_1\cap
E_*=\{0\}$ from that $e_0+\alpha e_0\in E_*$ implies $e_0=0$. Then,
$$e=P^{E_*}_{E_0}e+P^{E_0}_{E_*}e=(P^{E_*}_{E_0}e+\alpha(P^{E_*}_{E_0}e))+(P^{E_0}_{E_*}e-\alpha(P^{E_*}_{E_0}e)),\quad
\forall e\in E,$$

 which shows $E_1\oplus E_*\supset E$. So
$E_1\oplus E_*=E$. The proof ends.

\textbf{$6$ \quad Generalized regular point}

 Let $f$ be a $c^1$ map
from an open set $U$ in a Banach space $E$ into another Banach space.
  It is well known that immersion, submersion and subimmersion theorems hold for immersion, submersion and subimmersion points $x_0$ of $f$,respectively.
Each of them provides a vital approach with $f'(x_0)$ to study local behavior of $f$ near $X_0$.(Refer [Abr] and [Zei] .) The subimmersiom theorem is also said to be the rank theorem , which comes from its  assumption, $rank f'(x)=rank f'(x_0)< \infty $ for $x$ near $x_0.$ In [Ber] Berger shows that it is not yet known whether the rank theorem in advanced calculus holds even if  $f$ is Fredhlom map. In 1999 we gave the following  diffeomorphisms $u$ in $E$ and $v$ in $F$;$$ u(x)=T_0^+(f(x)-f(x_0))+(I_E-T_0^+T_0)(x-x_0)$$ and  $$v(y)=(f\circ u^{-1}\circ T_0^+)(y)+(I_F-T_0T_0^+)y$$, where $T_0=f'(x_0)$ and $T_0^+$ is a generalized inverse of $T_0$.(Refer [Ma $1$].) Then the following complete rank theorem in advanced calculus is established: suppose that $f'(x_0) (x_0\in U)$ is double splitting. There exist a neighborhood $U_0$ at $x_0$,$V_0$ at $o$,local diffeomorphisms $u:U_0\rightarrow u(U_0)$and $v:V_0\rightarrow v(V_0),$ such that $u(x_0)=0, u'(x_0)=I_E,v(x_0)=f(x_0), v'(0)=I_F$ and $$f(x)=(v\circ f'(x_0)\circ\circ u)(x) \quad \forall x\in U$$ if and only if $x_0$ is a generalized regular point of $f.$ (Refer [Ma$1$], Ma$9$], [Ber], and [Zei].) The theorem answers completely the problem presented by Berger,M.
In view of the above introduction and Theorem $3.2$ in this paper ones should recognize that the generalized regular point is a good mathematical concept.

\newpage
\vskip 0.1cm
\begin{center}{\bf References}
\end{center}
\vskip -0.1cm
\medskip
{\footnotesize
\def\REF#1{\par\hangindent\parindent\indent\llap{#1\enspace}\ignorespaces}

\REF{[Abr]}\ R. Abraham, J. E. Marsden, and T. Ratin, Manifolds,
Tensor Analysis and Applications, 2nd ed., Applied Mathematical
Sciences 75, Springer, New York, 1988.

\REF{[An]}\ V. I. Arnol'd, Geometrical Methods in the Theory of
Ordinary Differential Equations, 2nd ed., Grundlehren der
Mathematischen Wissenschaften 250, Springer, New York, 1988.

\REF{[Ber]}\ M, Brger, Nonlinearity and Functional Analysis, New York:
Academic Press, 1976.

\REF{[Bo]}B.Boss, D.D.Bleecker, Topology and Analysis, The
Atiyah-Singer Index Formula and Gauge-Theoretic Physics, New York
Springer-Verlag, 1985.

\REF{[Caf]}\ V. Cafagra, Global invertibility and finite
solvability, pp. 1-30 in Nonlinear Functional Analysis (Nework, NJ,
1987), edited by P. S. Milojevic, Lecture Notes in Pure and Appl.
Math. 121, Dekker, New York, 1990.

\REF{[H-M]}\ Qianglian Huang, Jipu Ma, Perturbation analysis of
generalized inverses of linear operators in Banach spaces, Linear
Algebra and its Appl., 389(2004), 359-364.

 \REF{[Ka]}\ Kato, T., Perturbation Theory for Linear
Operators, New York: Springer-Verlag, 1982.

\REF{[Ma1]}\ Jipu Ma, (1,2) inverses of operators between Banach
spaces and local conjugacy theorem, Chinese Ann. Math. Ser. B.
20:1(1999), 57-62.

\REF{[Ma2]}\ Jipu Ma, Rank theorem of operators between Banach
spaces, Sci China Ser. A 43(2000), 1-5.

\REF{[Ma3]}\ Jipu Ma, A generalized preimage theorem in global
analysis, Sci, China Ser. A 44:33(2001), 299-303.

\REF{[Ma4]}\ Jipu Ma, A rank theorem of operators between Banach
spaces, Front. of Math. in China 1(2006), 138-143.

\REF{[Ma5]} Jipu Ma, Three  classes of smooth Banach manifolds in
$B(E,F)$, Sci. China Ser. A 50(19)(2007), 1233-1239.

\REF{[Ma6]}\ Jipu Ma, A generalized transversility in global
analysis, Pacif. J. Math.,  236:2(2008),  357-371.

\REF{[Ma7]}\ Jipu Ma, A geometry characteristic of Banach spaces
with $C^1$-norm, Front. Math. China,  2014, 9(5): 1089-1103.

\REF{[Ma8]} Jipu Ma, A principle for critical point under
generalized regular constrain and Ill-posed Lagrange multipliers,
Numerical Functional and Optimization Theory 36.225:(2015), 369.

\REF{[Ma9]}\  Jipu Ma, Local conjugacy  theorem, rank theorem in
advanced calculus and a generalized principle for constructing Banach
manifold, Sci. China Ser. A 43 (2000), 1233-1237.

\REF{[N-C]}\ M. Z. Nashed, X. Chen, Convergence of Newton-like methods for singular equations using outer
 inverses, Numer. Math., 66:(1993), 235-257.

\REF{[P]}\ R. Penrose, A generalized inverse for Matrices, Proc. Cambridge
Philos. Soc. 1955, 51:406-413.

\REF{[Zei]}\ A. E. Zeilder, Nonlinear Functional Analysis and its Applications,
IV. Applications to Mathematical Physics. New York: Springer -verlag, 1988.

  1. Department of Mathematics, Nanjing University, Nanjing, 210093,
  P. R. China

  2. Tseng Yuanrong Functional Research Center, Harbin Normal
  University, Harbin, 150080, P. R. China

 E-mail address: jipuma@126.com

\end{document}